\documentclass[a4paper,12pt, table]{article}
\usepackage[text={170mm,254mm},centering]{geometry}

\usepackage[utf8]{inputenc}
\usepackage[T1]{fontenc}

\usepackage{amsmath} 
\usepackage{amsthm}
\usepackage{amsfonts}
\usepackage{mathtools}  
\usepackage{bm}

\usepackage{soul}

\usepackage{enumerate}
\usepackage{graphicx}
\usepackage{subcaption}
\usepackage{url}

\usepackage{tikz}
\usetikzlibrary{calc}
\usetikzlibrary{backgrounds}

\definecolor{myMagenta}{RGB}{255,0,255}
\definecolor{myCyan}{RGB}{0,255,255}
\tikzset{>=latex}

\usepackage{authblk}

\theoremstyle{plain}
\newtheorem{theorem}{Theorem}[section]

\newtheorem{corollary}[theorem]{Corollary}
\newtheorem{proposition}[theorem]{Proposition}
\newtheorem{conjecture}[theorem]{Conjecture}
\theoremstyle{definition}

\newtheorem{problem}{Problem}
\newtheorem{question}[problem]{Question}
\newtheorem{remark}{Remark}

\newcommand{\bZ}{\mathbb{Z}}

\usepackage[color links=true, backref=page]{hyperref}

\title{On a new $(21_4)$ polycyclic configuration}
\author[1]{Leah Wrenn Berman}
\author[2]{G\'{a}bor G\'{e}vay}
\author[3,4,5,6]{Toma\v{z} Pisanski}

\affil[1]{Department of Mathematics \& Statistics,
University of Alaska Fairbanks,
Fairbanks, USA}
\affil[2]{Bolyai Institute, University of Szeged, Szeged, Hungary}
\affil[3]{FAMNIT, University of Primorska, Koper, Slovenia}
\affil[4]{IAM, University of Primorska, Koper, Slovenia}
\affil[5]{Institute of Mathematics, Physics and Mechanics, Ljubljana, Slovenia}
\affil[6]{Faculty of Mathematics and Physics,
  University of Ljubljana, Ljubljana, Slovenia}

\date{\today}

\begin{document}

\maketitle



\begin{abstract}
When searching for small 4-configurations of points and lines, polycyclic configurations, in which every symmetry class of points and lines contains the same number of elements, 
have proved to be quite useful. In this paper we construct and prove the existence of a previously unknown $(21_4)$ configuration, which provides a counterexample to a conjecture 
of Branko Gr\"{u}nbaum. In addition, we study some of its most important properties; in particular, we make a comparison with the well-known Gr\"unbaum--Rigby configuration. 
We show that there are exactly two $(21_{4})$ geometric polycyclic   configurations and seventeen $(21_{4})$ combinatorial polycyclic configurations. We also discuss some 
possible generalizations.
\end{abstract}


\section{Introduction}

A breakthrough in the modern study of geometric configurations of points and lines came with the seminal paper~\cite{GruRig1990} 
of Gr\"{u}nbaum and Rigby in which the first geometric point-line representation of a 4-configuration was constructed. This $(21_4)$ 
configuration, which has 21 points and lines in which each point lies on 4 straight lines and each line passes through four points, was based on the work of Felix Klein \cite{Kl1878} on his famous quartic curve, and is nowadays known as the 
Gr\"{u}nbaum--Rigby configuration; we denote it by $\mathrm{GR}(21_{4})$. Later, Branko Gr\"{u}nbaum discovered a large 
number of $(n_4)$ configurations. Some of them were constructed in the spirit of $\mathrm{GR}(21_{4})$ (later called \emph{celestial} 
configurations), while others were constructed by various techniques from smaller ones. In 2003, Boben and Pisanski \cite{BoPi2003} 
initiated the theory of \emph{polycyclic configurations}, having $\mathrm{GR}(21_{4})$ and some other configurations from another 
paper of Gr\"{u}nbaum (co-authored by Harold Dorwart)~\cite{DoGr1992} as the prime models of such configurations.

A $(n_{k})$ combinatorial configuration is a collection of $n$ objects, called ``points'' and $n$ collections of ``points'', called  ``lines'', such that each point is incident with $k$ lines and each line contains $k$ points.  Each combinatorial configuration is in one-to-one correspondence to a bipartite graph in which each point and each line corresponds to a node of the configuration and incident point- and line-nodes are connected by an edge of the graph; this incidence graph is called the \emph{Levi graph} of the configuration.
If the points are distinct points in some Euclidean space (usually the plane) and the lines are distinct straight lines, then we call this a \emph{geometric} configuration, or a \emph{(strong) geometric realization} of the corresponding combinatorial configuration. A geometric configuration is \emph{polycyclic} if the orbits of the points and lines under the action of the maximal rotational symmetry group each have the same number of elements; we call these orbits the \emph{symmetry classes} of the configuration.

The study of polycyclic configurations 
was independently pursued and further developed by Gr\"{u}nbaum~\cite{Gru2009b} and Berman and coauthors (see, e.g.\ 
\cite{Ber2013,BeFaPi2016,BerJacVer2016,BeDeFaPiZi2020}). 
It is closely intertwined with graph theory as well; for details on this connection, see \cite{PisSer2013}.

\begin{figure}[htbp]
\begin{center}
\subcaptionbox{
}[.4\textwidth]{
\begin{tikzpicture}[vtx/.style={draw, circle, inner sep = 1.5 pt}, lbl/.style={midway, inner sep = 1 pt, fill =white, font = \footnotesize}, lin/.style={draw, square, inner sep = 12 pt}] 
\pgfmathsetmacro{\rr}{3.5}
\pgfmathsetmacro{\rad}{\rr*cos(2*180/7)/cos(180/7)}
\pgfmathsetmacro{\radd}{\rr*cos(2*180/7)/cos(180/7)*cos(3*180/7)/cos(2*180/7)}
\foreach \i in {0,..., 6}{
\node[vtx, fill = red](u\i) at (2*\i*180/7:\rr){};
\node[vtx, fill=blue] (v\i) at (2*\i*180/7+180/7:\rad){};
\node[vtx, fill=green!60!black] (w\i) at (2*\i*180/7:\radd){};
}
\begin{scope}[on background layer]
\foreach \i in {0,...,6}{
\draw[green!60!black, thick] let \n1 = {int(mod(\i+2,7))} in (u\i) -- (u\n1); 
\draw[red, thick] let \n1 = {int(mod(\i+3,7))} in (v\i) -- (v\n1); 
\draw[blue, thick] let \n1 = {int(mod(\i+3,7))} in (u\i) -- (u\n1); 
}
\end{scope}

\end{tikzpicture}
}
\hspace{1cm}
\subcaptionbox{
}[.4\textwidth]{
\begin{tikzpicture}[vtx/.style={draw, circle, inner sep = 2 pt}, lbl/.style={midway, inner sep = 1 pt, fill =white, font = \footnotesize}, lin/.style={draw, inner sep = 2.3 pt}] 
\pgfmathsetmacro{\rr}{3.5}
\foreach \i in {0,..., 6}{
\node[vtx, fill = red] (R\i) at (2*\i*180/7:{\rr}){};
\node[lin, fill = green!60!black] (r\i) at (2*\i*180/7+180*2/7:{\rr*4.5/6}){}; 
\node[vtx, fill = blue] (B\i) at (2*\i*180/7+180/7:{\rr*4/6}){};
\node[lin, fill = red] (b\i) at (2*\i*180/7 +4*180/7:{\rr*3.5/6}){}; 
\node[vtx, fill = green!60!black] (G\i) at (2*\i*180/7 + 2*180/7:{\rr*2/6}){};
\node[lin, fill = blue] (g\i) at (2*\i*180/7+3*180/7:{\rr*1/6}){}; 
}
\foreach \i in {0,...,6}{
\draw (R\i) -- (r\i);
\draw let \n1 = {int(mod(\i+2,7))} in (r\i) -- (R\n1);
\draw (r\i) -- (B\i);
\draw let \n1 = {int(mod(\i+1,7))} in (r\i) -- (B\n1);

\draw (b\i) -- (B\i);
\draw let \n1 = {int(mod(\i+3,7))} in (b\i) -- (B\n1);
\draw (b\i) -- (G\i);
\draw let \n1 = {int(mod(\i+2,7))} in (b\i) -- (G\n1);

\draw (g\i) -- (G\i);
\draw let \n1 = {int(mod(\i+1,7))} in (g\i) -- (G\n1);
\draw (g\i) -- (R\i);
\draw let \n1 = {int(mod(\i+3,7))} in (g\i) -- (R\n1);

}

\end{tikzpicture}
}
\caption{(a) The Gr\"{u}nbaum--Rigby $(21_4)$ geometric configuration, denoted by $\mathrm{GR}(21_4)$, and (b) its Levi graph. In the Levi graph, point-vertices are shown with 
circles and line-vertices are shown with squares, and the colours correspond to the colours of the symmetry classes of the points and lines of the configuration.
We note that the colours are also consistent with the self-dualities (in the sense that self-duality maps of the configuration preserve the colours; these maps correspond to certain 
automorphisms of the Levi graph).}
\label{fig:GR}
\end{center}
\end{figure}
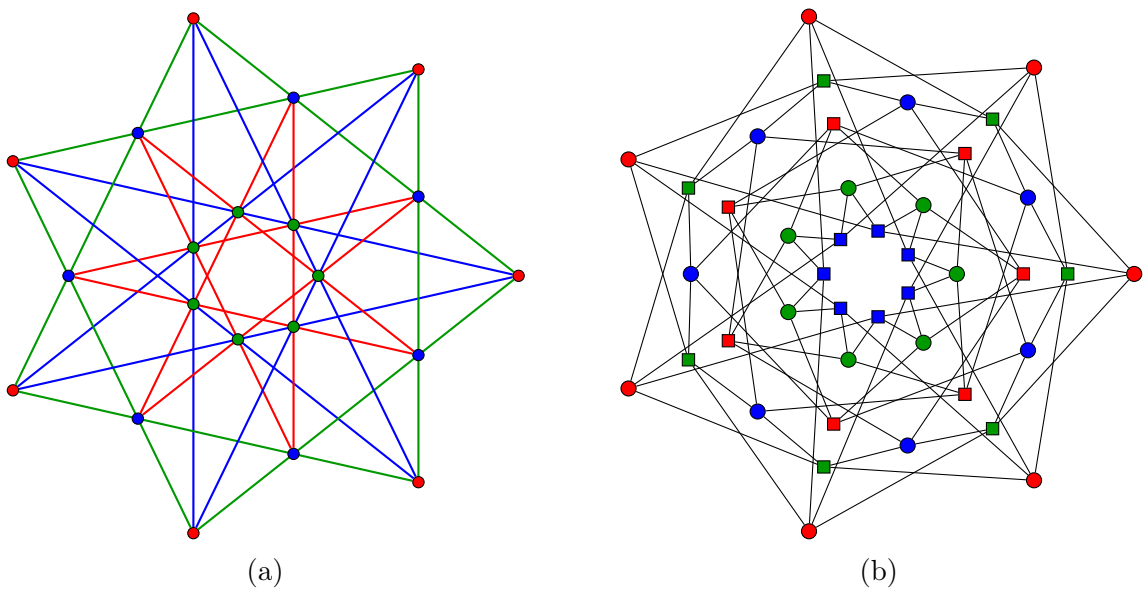

In their paper, Gr\"{u}nbaum and Rigby conjectured:
\begin{enumerate}
\item that no other $(21_4)$ configuration exists, and
\item no $(n_{4})$ configuration exists for  $n < 21$. 
\end{enumerate}
It was a big surprise 
when Gr\"{u}nbaum himself disproved the second part of this conjecture \cite{Gru2008a} by constructing a $(20_4)$ configuration, which we denote by 
$\mathrm{G}(20_{4})$. At that time, it was widely believed that the $\mathrm{GR}(21_{4})$ configuration was the only geometric $4$-configuration for 
$n = 21$, and that $\mathrm{G}(20_{4})$ is the smallest geometric 4-configuration. However, a few years later, J\"urgen Bokowski and his coauthors 
showed that there are no $(n_{4})$ configurations for $n\leq 17$~\cite{BokSch2005}, that there are exactly two distinct $(18_4)$ configurations~\cite
{BokSch2013,BokPil2014}, and that no geometric $(19_4)$ configuration exists~\cite{BokPil2015}.

\begin{figure}[h]
\begin{center}
\includegraphics[width=0.6\textwidth, angle = -90]{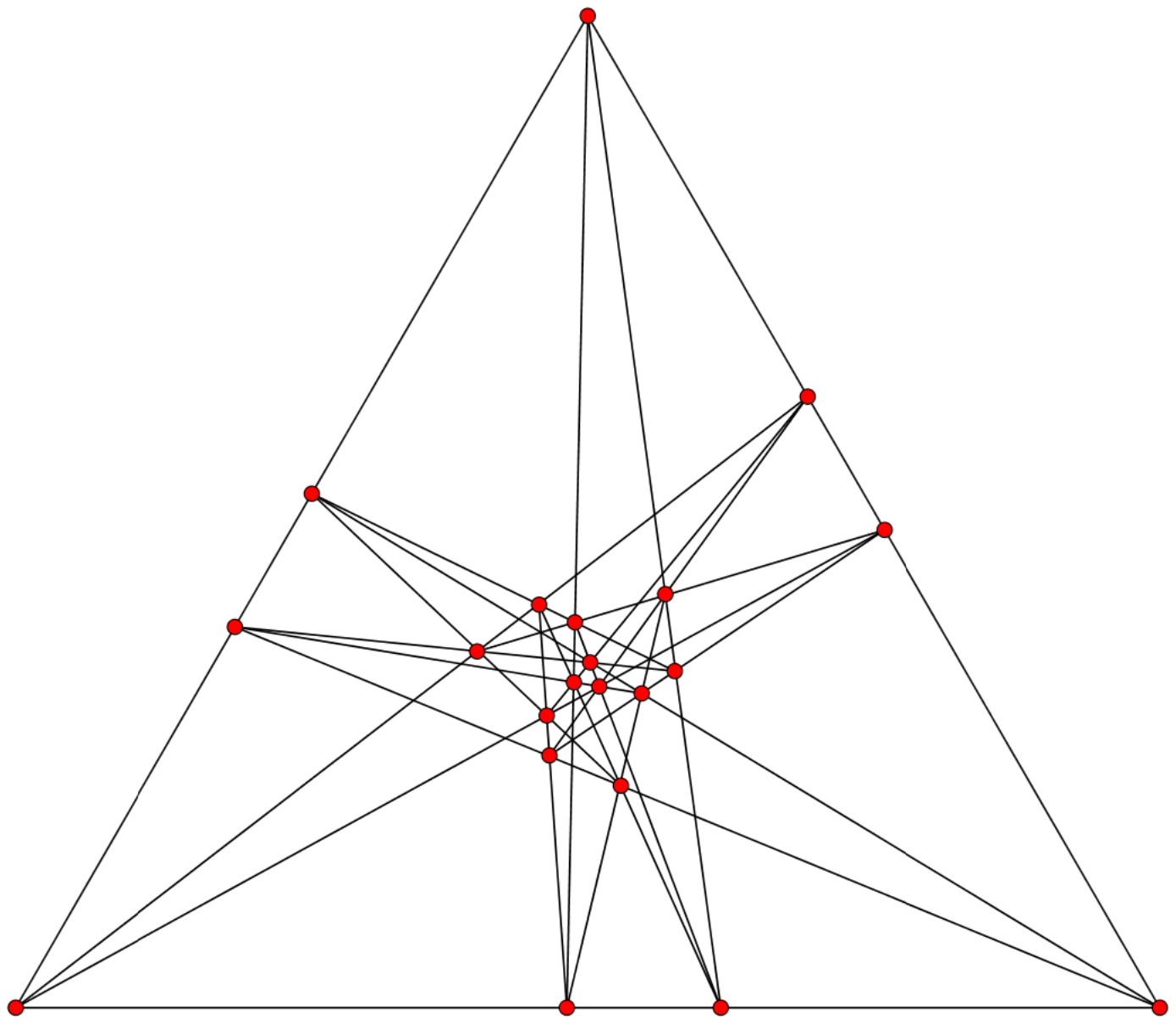}
\caption{A new $(21_4)$ geometric configuration, denoted by $\mathrm{B}(21_{4})$.}
\label{fig:B21}
\end{center}
\end{figure}

A number of months ago, the first author of this paper constructed a new $(21_4)$ geometric configuration, depicted in Figure~\ref{fig:B21}, which provides a counterexample 
to the first part of the Gr\"{u}nbaum--Rigby conjecture. We denote this configuration by $\mathrm{B}(21_{4})$. The main goal of this paper is to provide a proof of existence of 
this configuration; we present both a synthetic and an analytic proof, since we believe that both have their benefits, and may form a suitable basis for extending the research to 
configurations with analogous structure. For the same reason,  we also discuss some interesting structural properties of this configuration.


\section{A comparison of the configurations $\mathrm{GR}(21_{4})$ and $\mathrm{B}(21_4)$} \label{sect:comp}


It is not hard to verify that $\mathrm{B}(21_{4})$ is  combinatorially distinct from $\mathrm{GR}(21_{4})$.  Namely, one can compute the Levi graphs (that is, the point-line incidence 
graphs, with one vertex of the graph for each point and line of the configuration, with a point-vertex connected by an edge to a line-vertex if and only if the point and line are incident 
in the configuration) of both configurations. We used the computer algebra system \emph{Sage} to prove that the two 4-valent graphs on 42 vertices are non-isomorphic. For instance, 
the Levi graph of $\mathrm{B}(21_{4})$ has only 12 automorphisms, while the Levi graph of $\mathrm{GR}(21_{4})$ has 672 automorphisms (including bipartition-reversing automorphisms, 
which correspond to self-dualities of these configurations).

The Levi graph of the Gr\"{u}nbaum--Rigby configuration, which we denote by L(GR), can be described as the Kronecker cover over the line graph of the renowned Heawood graph. 
Its automorphism group contains 672 elements. Half of them correspond to combinatorial self-dualities, while the other half correspond to combinatorial automorphisms of 
$\mathrm{GR}(21_{4})$. As for the latter, we know that the automorphism group of both the Heawood graph and the Gr\"unbaum--Rigby 
configuration is $\mathbf{PGL}(2,7)$ of order 336~\cite{GruRig1990}, and is a subgroup of index 2 in the automorphism group of L(GR). 
We observe that out of the 336 combinatorial symmetries, only 14 are geometrically realizable in the standard polycyclic realization; also, from the 336 combinatorial self-dualities, 14 are geometrically realizable. This means that $\mathrm{GR}(21_{4})$ shown in  Figure~\ref{fig:GR} geometrically realizes 28 out of the 672 graph automorphisms. 

We used programs written in \emph{Sage} 
to compute all semi-regular automorphisms of L(GR) and the corresponding quotient graphs. The quotients that are bipartite correspond to \emph{reduced Levi graphs}. We often 
abbreviate reduced Levi graphs as RLG. For their definition and some properties, including voltage groups and voltage assignments, see~\cite{Gru2009b, PisSer2013, Ber2013}). 
Our computations show that there are 314 semi-regular automorphisms producing 8 distinct quotient graphs of L(GR). However, only two of them are bipartite, hence there exist 
only two non-isomorphic RLGs that can possibly correspond to polycyclic realizations  of of the Gr\"unbaum-Rigby configuration. 

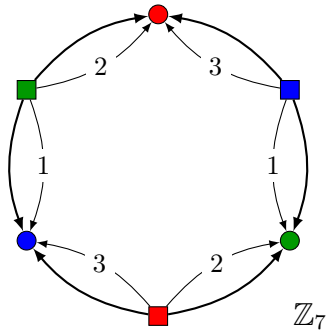
\begin{figure}[h]
\begin{center}
\begin{tikzpicture}[vtx/.style={draw, circle, inner sep = 2.5 pt}, lbl/.style={midway, inner sep = 3 pt, fill =white, font = \footnotesize}, lin/.style={draw, inner sep = 3.5 pt}] 
\node[vtx, fill = red] (u) at (90:2){};
\node[lin, fill = green!60!black] (L) at (90+60:2){};  
\node[vtx, fill = blue] (v) at (90+120:2){};
\node[lin, fill = red] (M) at (90+120+60:2){}; 
\node[vtx, fill = green!60!black] (w) at (90+240:2){};
\node[lin, fill = blue] (N) at (90+240+60:2){}; 

\draw[<-] (u) to[bend left = 20] node[lbl]{2} (L);
\draw[<-,  thick] (u) to[bend right = 20]  (L);
\draw[<-] (v) to[bend right = 20] node[lbl]{1} (L);
\draw[<-,  thick] (v) to[bend left = 20]  (L);

\draw[<-] (v) to[bend left = 20] node[lbl]{3} (M);
\draw[<-,  thick] (v) to[bend right = 20]  (M);
\draw[<-] (w) to[bend right = 20] node[lbl]{2} (M);
\draw[<-,  thick] (w) to[bend left = 20]  (M);

\draw[<-] (w) to[bend left = 20] node[lbl]{1} (N);
\draw[<-,  thick] (w) to[bend right = 20]  (N);
\draw[<-] (u) to[bend right = 20] node[lbl]{3} (N);
\draw[<-,  thick] (u) to[bend left = 20]  (N);

\node at ($(M)+(2,0)$) {$\mathbb{Z}_{7}$};

\end{tikzpicture}

\caption{The reduced Levi graph RLG(GR) with voltage group $\bZ_7$ for the polycyclic Gr\"{u}nbaum--Rigby 
$\mathrm{GR}(21_{4})$ configuration with seven-fold rotational symmetry. The colours for the symmetry classes match the colours from Figure \ref{fig:GR}.}
\label{fig:sRLG}
\end{center}
\end{figure}

The first RLG, on 6 vertices, is expected, since it can be deduced from Figure~\ref{fig:GR}
, and is depicted in Figure~\ref{fig:sRLG}. The associated voltage group is $\bZ_7$, 
consistent with the 7-fold rotational symmetry of the geometric Gr\"unbaum--Rigby configuration. 

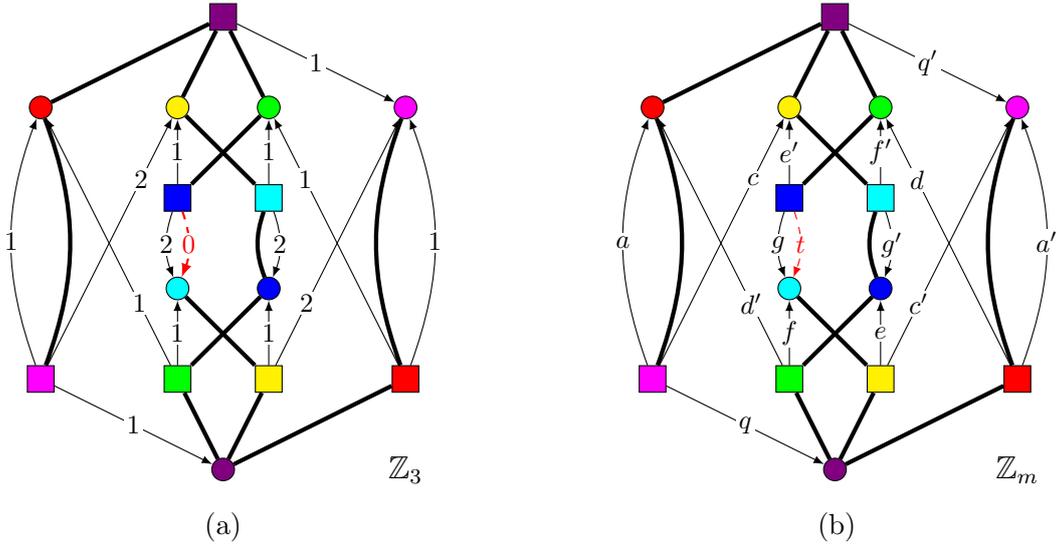
\begin{figure}[!h]
\begin{center}
\subcaptionbox{
\label{fig:RLG(B)}}[.4\textwidth]{
\begin{tikzpicture}[scale = 1.2, lbl/.style={midway, inner sep = 1 pt, fill =white, font = \footnotesize}]
\node[draw, , fill = violet, inner sep = 5 pt] (blackLin) at (0,0) {}; 
\node[draw, circle, fill = red, inner sep = 3 pt] (redPt) at (-2, -1) {};
\node[draw, circle, fill = yellow, inner sep = 3 pt, ] (yellowPt) at (-.5, -1) {};
\node[draw, circle, fill = green, inner sep = 3 pt, ] (greenPt) at (.5, -1) {};
\node[draw, circle, fill = myMagenta, inner sep = 3 pt, ] (magentaPt) at (2, -1) {};
\node[draw, , fill = blue, inner sep = 5 pt, ] (cyanLin) at (-.5, -2) {}; 
\node[draw, , fill = myCyan, inner sep = 5 pt, ] (blueLin) at (.5, -2) {}; 
\node[draw, circle, fill = myCyan, inner sep = 3 pt, ] (cyanPt) at (-.5, -3) {};
\node[draw, circle, fill = blue, inner sep = 3 pt, ] (bluePt) at (.5, -3) {};
\node[draw, , fill = myMagenta, inner sep = 5 pt, ] (redLin) at (-2, -4) {}; 
\node[draw, , fill = red, inner sep = 5 pt, ] (magentaLin) at (2, -4) {}; 
\node[draw, , fill = green, inner sep = 5 pt, ] (greenLin) at (-.5, -4) {}; 
\node[draw, , fill = yellow, inner sep = 5 pt, ] (grayLin) at (.5, -4) {}; 
\node[draw, circle, fill = violet, inner sep = 3 pt] (orangePt) at (0,-5) {}; 

\foreach \i/\j in {blackLin/redPt,blackLin/yellowPt, blackLin/greenPt,yellowPt/blueLin, greenPt/cyanLin, cyanPt/grayLin, bluePt/greenLin, orangePt/magentaLin}{
\draw[ultra thick] (\i) -- (\j);}
\draw[ultra thick] (redLin) to[bend right = 20] (redPt);
\draw[ultra thick] (magentaLin) to[bend left = 20] (magentaPt);
\draw[thick, dashed, red, ->] (cyanLin) to[bend left = 20]node[lbl]{$0$}  (cyanPt); 
\draw[ultra thick] (blueLin) to[bend right = 20] (bluePt);
\draw[ultra thick] (orangePt) to[]  (greenLin);
\draw[ultra thick] (orangePt) to[] (grayLin);

\draw[<-] (redPt) to[bend right=20] node[lbl]{$1$} (redLin); 
\draw[<-] (magentaPt) to[bend left=20] node[lbl]{$1$} (magentaLin); 
\draw[<-] (yellowPt) to[] node[lbl, near start]{$2$} (redLin); 
\draw[<-] (greenPt) to[] node[lbl, near start]{$1$} (magentaLin);
\draw[<-] (yellowPt) to[] node[lbl]{$1$} (cyanLin);
\draw[<-] (greenPt) to[] node[lbl]{$1$} (blueLin);
\draw[<-] (cyanPt) to[bend left=20] node[lbl]{$2$} (cyanLin);
\draw[<-] (bluePt) to[bend right=20] node[lbl]{$2$} (blueLin);
\draw[<-] (cyanPt) to[] node[lbl]{$1$} (greenLin);
\draw[<-] (bluePt) to[] node[lbl]{$1$} (grayLin);
\draw[<-] (redPt) to[] node[lbl, near end]{$1$} (greenLin);
\draw[<-] (magentaPt) to[] node[lbl, near end]{$2$} (grayLin);

\draw[<-] (orangePt) to[] node[lbl]{$1$} (redLin); 
\draw[->] (blackLin) to[] node[lbl]{$1$} (magentaPt); 

\draw (2,-5) node {$\mathbb{Z}_{3}$};

\end{tikzpicture}
}
\hspace{1cm}
\subcaptionbox{
\label{fig:RLG-general}}[.4\textwidth]{
\begin{tikzpicture}[scale = 1.2, lbl/.style={midway, inner sep = 1 pt, fill =white, font = \footnotesize}]
\node[draw, , fill = violet, inner sep = 5 pt] (blackLin) at (0,0) {}; 
\node[draw, circle, fill = red, inner sep = 3 pt] (redPt) at (-2, -1) {};
\node[draw, circle, fill = yellow, inner sep = 3 pt, ] (yellowPt) at (-.5, -1) {};
\node[draw, circle, fill = green, inner sep = 3 pt, ] (greenPt) at (.5, -1) {};
\node[draw, circle, fill = myMagenta, inner sep = 3 pt, ] (magentaPt) at (2, -1) {};
\node[draw, , fill = blue, inner sep = 5 pt, ] (cyanLin) at (-.5, -2) {}; 
\node[draw, , fill = myCyan, inner sep = 5 pt, ] (blueLin) at (.5, -2) {}; 
\node[draw, circle, fill = myCyan, inner sep = 3 pt, ] (cyanPt) at (-.5, -3) {};
\node[draw, circle, fill = blue, inner sep = 3 pt, ] (bluePt) at (.5, -3) {};
\node[draw, , fill = myMagenta, inner sep = 5 pt, ] (redLin) at (-2, -4) {}; 
\node[draw, , fill = red, inner sep = 5 pt, ] (magentaLin) at (2, -4) {}; 
\node[draw, , fill = green, inner sep = 5 pt, ] (greenLin) at (-.5, -4) {}; 
\node[draw, , fill = yellow, inner sep = 5 pt, ] (grayLin) at (.5, -4) {}; 
\node[draw, circle, fill = violet, inner sep = 3 pt] (orangePt) at (0,-5) {}; 

\foreach \i/\j in {blackLin/redPt,blackLin/yellowPt, blackLin/greenPt,yellowPt/blueLin, greenPt/cyanLin, cyanPt/grayLin, bluePt/greenLin, orangePt/magentaLin}{
\draw[ultra thick] (\i) -- (\j);}
\draw[ultra thick] (redLin) to[bend right = 20] (redPt);
\draw[ultra thick] (magentaLin) to[bend left = 20] (magentaPt);
\draw[red, dashed, ->] (cyanLin) to[bend left = 20] node[lbl]{$t$} (cyanPt);
\draw[ultra thick] (blueLin) to[bend right = 20] (bluePt);
\draw[ultra thick] (orangePt) to[]  (greenLin);
\draw[ultra thick] (orangePt) to[] (grayLin);

\draw[<-] (redPt) to[bend right=20] node[lbl]{$a$} (redLin);
\draw[<-] (magentaPt) to[bend left=20] node[lbl]{$a'$} (magentaLin);
\draw[<-] (yellowPt) to[] node[lbl, near start]{$c$} (redLin);
\draw[<-] (greenPt) to[] node[lbl, near start]{$d$} (magentaLin);
\draw[<-] (yellowPt) to[] node[lbl]{$e'$} (cyanLin);
\draw[<-] (greenPt) to[] node[lbl]{$f'$} (blueLin);
\draw[<-] (cyanPt) to[bend left=20] node[lbl]{$g$} (cyanLin);
\draw[<-] (bluePt) to[bend right=20] node[lbl]{$g'$} (blueLin);
\draw[<-] (cyanPt) to[] node[lbl]{$f$} (greenLin);
\draw[<-] (bluePt) to[] node[lbl]{$e$} (grayLin);
\draw[<-] (redPt) to[] node[lbl, near end]{$d'$} (greenLin);
\draw[<-] (magentaPt) to[] node[lbl, near end]{$c'$} (grayLin);

\draw[<-] (orangePt) to[] node[lbl]{$q$} (redLin);
\draw[<-] (magentaPt) to[] node[lbl]{$q'$} (blackLin);

\draw (2,-5) node {$\mathbb{Z}_{m}$};

\end{tikzpicture}
}
\caption{(a) The reduced Levi graph RLG(B) for the configuration $\mathrm{B}(21_{4})$. The voltage group is 
$\mathbb{Z}_{3}$, expressing the fact that this configuration exhibits threefold rotational symmetry. Point 
orbits of $\mathrm{B}(21_{4})$ are represented by circular nodes and line orbits by rectangular nodes. 
The colour-preserving half-turn symmetry of the graph corresponds to a self-duality of the configuration. (b) A version of the RLG with generic parameters. }
\end{center}
\end{figure}

However, the second one, depicted in Figure~\ref{fig:RLG(B)}, was quite unexpected. It has 14 vertices, and the corresponding voltage group 
is $\bZ_3$. Initially, we wanted to know if there existed a polycyclic geometric realization of the Gr\"{u}nbaum--Rigby configuration with 3-fold 
rotational symmetry. All our attempts to generate such a realization based on the RLG shown in Figure~\ref{fig:RLG(B)} failed (see Section \ref{sect:gen}).

There is a simple algorithm that can produce a Levi graph from a reduced Levi graph. For a reduced Levi graph with voltage group $\mathbb{Z}_{m}$, the notation 
\tikz{\node[draw, inner sep = 1.5 pt](L) at (0,0){$L$}; \node[draw, circle,inner sep = 1 pt](v) at (1.5,0){$v$}; \draw[->] (L) --node[midway, inner sep = 1 pt, fill= white]{$a$} (v)} 
means that there is a symmetry class of points labeled $v$, with elements $v_{i}$, $i = 0,\ldots, m-1$; a symmetry class of lines $L$, with elements $L_{i}$, $i = 0,\ldots,m-1$; 
and that each line $L_{i}$ is incident with vertex $v_{i+a}$, with index arithmetic taken modulo $m$. (A more detailed description of the relationship between reduced Levi graphs 
with voltage group $\mathbb{Z}_{m}$ can be found in \cite{BoPi2003, PisSer2013, Ber2013}). For convenience, we provide an incidence table from the reduced Levi graph shown 
in Figure \ref{fig:RLG(B)}, see Table~\ref{tab:inci}.

\begin{table}[!h]
\begin{center}
$$\begin{array}{||c|cccc||c|cccc||c|cccc||}
r_0  & M_0 & M_1 & G_1 & P_0 & r_1  & M_1 & M_2 & G_2 & P_1 & r_2 & M_2 & M_0  & G_0 & P_2\\
y_0  & M_2 & B_1 & C_0 & P_0  & y_1 & M_0 & B_2 & C_1 & P_1 & y_2 & M_1 & B_0  & C_2 & P_2 \\
g_0  & B_0 & C_1 & R_1 & P_0 & g_1  & B_1 & C_2 & R_2 & P_1 & g_2  & B_2 & C_0  & R_0 & P_2 \\
m_0 & Y_2 & R_0 & R_1 & P_1 & m_1 & Y_0 & R_1 & R_2 & P_2 & m_2 & Y_1 & R_2  & R_0 & P_0 \\
b_0  & C_0 & C_2 & Y_1 & G_0 & b_1 & C_1 & C_0 & Y_2 & G_1 & b_2  & C_2 & C_1  & Y_0 & G_2 \\
c_0  & B_0 & B_2 & G_1 & Y_0 & c_1 & B_1 & B_0 & G_2 & Y_1  & c_2  & B_2 & B_1 & G_0 & Y_2 \\
p_0 & R_0 & Y_0 & G_0 & M_1 & p_1 & R_1 & Y_1 & G_1 & M_2 & p_2  & R_2 & Y_2 & G_2 & M_0 \\
\end{array}
$$
\caption{The incidence table of $\mathrm B(21_4)$.}
\label{tab:inci}
\end{center}
\end{table}


\section{A synthetic proof of the existence of $\mathrm B(21_4)$} \label{sect:synthetic}

\bigskip

\begin{theorem} \label{thm:MainThm}
There exists a self-dual geometric, polycyclic $(21_4)$ configuration with threefold rotational symmetry.
\end{theorem}

Before proving this theorem, we need some preliminaries.

\subsection{Quasi-configurations}

Recall that an \emph{incidence structure} $\mathcal C$ is a triple $\mathcal C = (P, B, I )$, where $P$ is the set of \emph{points}, \emph{B} is the 
set of \emph{lines} (or \emph{blocks}), and $I \subseteq P \times B$ is the \emph{incidence relation}~\cite{PisSer2013}. Given an incidence stucture 
$\mathcal{C}= (P, B, I )$, assume that $P$ is a disjoint union of subsets $P_i$ of cardinality $p^i$ $(i =1 ,2, \dots, m(P))$, and $L$ is a disjoint union 
of subsets $L_j$ of cardinality $n^j$ $(j =1 ,2, \dots, m(L))$. We call $\mathcal{C}$ a (combinatorial) \emph{quasi-configuration of type} 
$$
\bigl(\bigl(p^1_{q_1}\bigr)\bigl(p^2_{q_2}\bigr)\dots\bigl(p^{m(P)}_{q_{m(P)}}\bigr),
\bigl(n^1_{k_1}\bigr)\bigl(n^2_{k_2})\dots\bigl(n^{m(L)}_{k_{m(L)}}\bigr)\bigr),
$$
if each point in $P_i$ is incident with $q_i$ lines, and each line in $L_j$ is incident with $n_j$ points. (We note that here we adopt the term used by Bokowski 
and Pilaud~\cite{BokPil2016}, but with slightly different meaning.) 

Observe that a quasi-con\-fig\-uration with $m(P)=m(L)=1$ is a configuration in the usual sense. In analogy with the case of configurations, if all the numerical 
parameters of a quasi-configuration $\mathcal C$ are the same for the points and the lines, we use the simplified notation
$$
\bigl(\bigl(n^1_{k_1}\bigr)\bigl(n^2_{k_2})\dots\bigl(n^{m(L)}_{k_{m(L)}}\bigr)\bigr),
$$
and we say that $\mathcal C$ is \emph{balanced} (here we adopt the term introduced by Gr\"unbaum~\cite{Gru2009b}). 
In particular, below we construct a quasi-configuration of type $((6_2)(9_4))$, where the type notation shows that it contains 
\begin{itemize}
\item 6 points, each incident to 2 lines, 
\item 9 points, each incident to 4 lines, 
and conversely, 
\item 6 lines, each incident to 2 points and 
\item 9 lines, each incident to 4 points.
\end{itemize}

\subsection{Self-reciprocity}

\begin{figure}[!h]
\begin{center}
\includegraphics[width=0.4\textwidth]{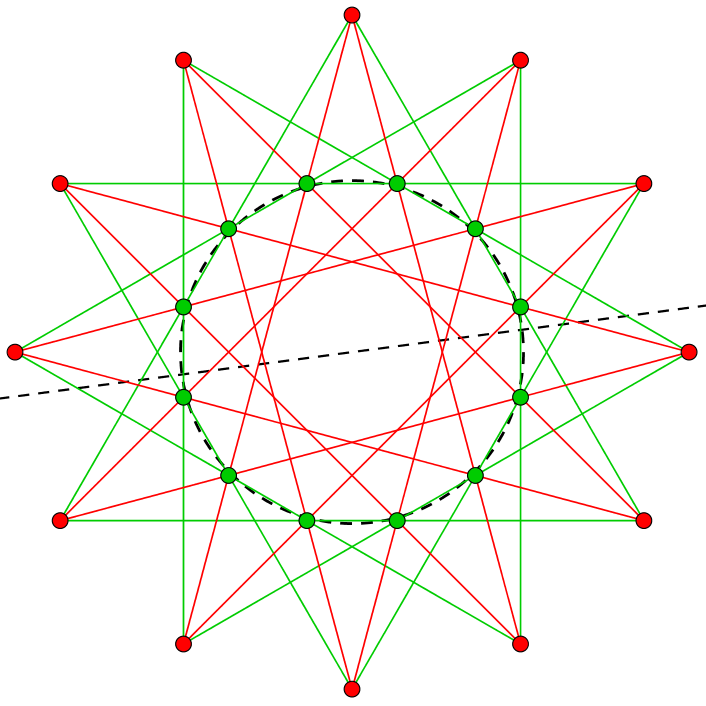}
\caption{Example of a configuration which is both rotationally and reflexibly self-reciprocal. The circle of reciprocity is shown (drawn dashed). 
The angle of rotation for attaining coincidence is $\pi/12$. The reciprocal image will also coincide with the original copy by reflecting it in the 
mirror line shown dashed (there are 12 such mirror lines).}
\label{fig:RotRef}
\end{center}
\end{figure}

It is clear that the notions of duality and self-duality of configurations applies also to quasi-configurations. In particular, it is also clear that a self-dual 
quasi-configuration is necessarily balanced. A stronger version of duality is when it is induced by reciprocity with respect to a circle. We say that a 
quasi-configuration $\mathcal C$ is \emph{self-reciprocal} if there is a circle $\Omega$ such that the reciprocity with respect to $\Omega$ sends 
$\mathcal C$ to its isometric copy. We distinguish three particular cases of self-reciprocity. We call $\mathcal C$ \emph{perfectly self-reciprocal} if 
it coincides with its reciprocal. The Gr\"unbaum--Rigby configuration provides an example of a perfectly self-reciprocal configuration (further examples
occur in~\cite{GevPok2023}, where this notion is applied for the first time).  Two slightly weaker versions are when for attaining coincidence one has 
to apply a subsequent rotation or reflection on the reciprocal image. In these cases we speak of a \emph{rotationally} or \emph{reflexibly self-reciprocal} 
configuration, respectively. Figure~\ref{fig:RotRef} shows an example of a polycyclic $(24_4)$ con\-figuration which is rotationally self-reciprocal. One 
directly observes that, in  addition, it is mirror symmetric (with 12 mirror lines); this implies that it is reflexibly self-reciprocal as well. We remark that 
Branko Gr\"unbaum in his book~\cite{Gru2009b} uses the term \emph{oppositely selfpolar} for a configuration with this latter property; in his Figure 
5.8.2 he presents precisely this configuration for an example.
 
Given a polycyclic (quasi-)configuration $\mathcal C$ possessing any of the self-reciprocity properties mentioned above, any 
(say, the $i$th) orbit of points of  $\mathcal C$ is located on a circle $\Gamma_i$ concentric with $\Omega$, and the orbit of the 
corresponding polar lines has an incircle $\Gamma'_i$ (also concentric with $\Omega$). Obviously, $\Gamma_i$ and $\Gamma'_i$ 
are inverse images of each other with respect to $\Omega$, in other words, $\Omega$ is the \emph{midcircle} of these circles. 

The mid-circle can be constructed in the following way. Take a ray starting from the common centre $O$, and let it intersect 
$\Gamma_i$ and $\Gamma'_i$ in points $P$ and $P'$, respectively. Take the Thales' circle $\mathcal T$ with diameter 
$PP'$, and construct a tangent line to $\mathcal T$ from the point $O$. Let $T$ be the point of tangency. Then the midcircle 
is obtained as a circle of radius $OT$ centred at $O$. We recall that this construction is based on some elementary properties 
of the inversion~\cite{CoxGre}, which can be summarized in the following proposition. 

\begin{proposition}
Let $\Omega$ be the circle of an inversion $\varphi$, and let $\Gamma_1$ be a circle. Then the following conditions are equivalent. 
\begin{enumerate}
\item 
$\Gamma_1$ is invariant under $\varphi$;
\item
if $\Gamma_1$ passes through a point $P$, then it also passes through $\varphi(P)$;
\item
$\Gamma_1$ is orthogonal to $\Omega$.
\end{enumerate}
\end{proposition}

\subsection{A quasi-configuration of type $((6_2)(9_4))$}

\begin{proposition} \label{prop:QCB}
There exists a self-reciprocal quasi-configuration of type $((6_2)(9_4))$ with threefold rotational symmetry. It is movable, with one degree of freedom.
\end{proposition}

\begin{figure}[!h] 
\begin{center}
\includegraphics[width=0.825\textwidth]{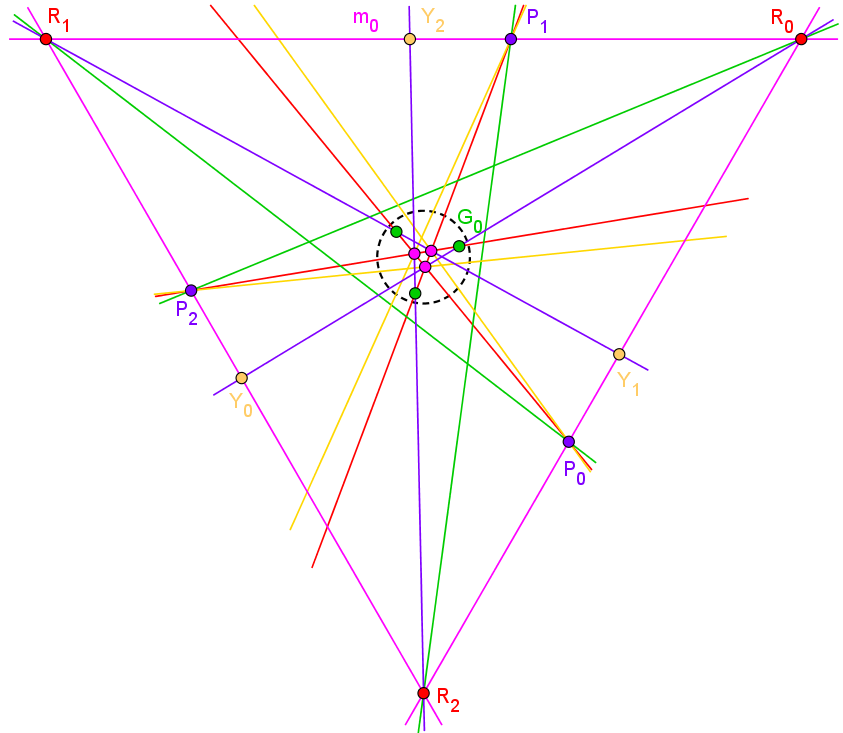}
\caption{The quasi-configuration QC(B). The circle of reciprocity is also shown (drawn dashed).}
\label{fig:quasi}
\end{center}
\end{figure}

We denote this quasi-configuration by QC(B). It is depicted in Figure~\ref{fig:quasi}. We distinguish the orbits of points and lines (w.r.t.\ its rotational 
symmetry group) by colours, namely, we use green, magenta, purple, red and yellow. The points and lines will be denoted accordingly by $X_i$ and 
$x_i$ $(i = 0, 1, 2)$, respectively, where $X$ and $x$ is the initial of the name of the corresponding colour. The corresponding orbit of points and lines
will be denoted by $(X)$ and $(x)$, respectively. As we shall see below, this quasi-configuration forms a substructure of $\mathrm{B}(21_4)$; hence 
here we use for labelling its points and lines the labels taken from Table \ref{tab:inci}.

\begin{proof}
We give a construction for QC(B) in the following 11 steps. 
\begin{enumerate} [(1)]
\item \label{step:initialstep}
Fix an equilateral triangle with red vertices $R_i$ and magenta side lines $m_i=R_{i}R_{i+1}$. 
\end{enumerate}

\begin{remark}
Here, and throughout the construction, the indices are meant modulo 3. In addition, we use the convention that rotation either of a point $X_i$ or of a line $x_i$ by angle 
$+120^{\circ}$ (i.e., counterclockwise) increases $i$ by one; the centre of rotation coincides with centre $O$ of the triangle $R_0R_1R_2$.
\end{remark}

\begin{enumerate} [(1)]
\addtocounter{enumi}{1}
\item \label{step:yellowpoint} 
Take a yellow point $Y_2$ on the line $m_0$ such that it can be shifted freely in the interior of the segment $R_1M_{01}$, where $M_{01}$ denotes the midpoint of the side $R_0R_1$
of the triangle $R_0R_1R_2$. Take also the rotates of this point by angle $\pm 120^{\circ}$.
\item 
Take the purple line $p_2\coloneqq Y_2R_2$ and its corresponding rotates.
\item \label{step:aux} 
Let $A$ be an auxiliary point defined as the intersection $A\coloneqq  \overline{Y_0Y_1} \cap \mathcal C(R_1OR_2)$, where $\overline{Y_0Y_1}$ is the segment connecting
$Y_0$ and $Y_1$, and $\mathcal C(R_1OR_2)$ is an auxiliary circle determined by $R_1$, $R_2$ and the centre $O$.
\item \label{step:greenlines} 
Take the green line $g_1\coloneqq  R_2A$ and its corresponding rotates.
\item \label{step:purplepoints}
Take the purple point of intersection $P_0\coloneqq g_0 \cap m_2$, and its rotates $P_1$ and $P_2$.
\item \label{step:mid} 
Take the circumcircle of the points $P_0$, $P_1$, $P_2$, and the incircle of the triangle formed by the lines $p_0$, $p_1$, $p_2$.
Construct the mid-circle of these circles; denote it by $\Omega$. 
\item \label{step:redlines} 
Take the circumcircle $\gamma_R$ of the points $R_0$, $R_1$, $R_2$, and invert it in the circle  $\Omega$. 
Denote the inverse circle by $\gamma_r$.
Draw tangents to $\gamma_r$ from the point $P_1$, and choose the one subtending the smaller angle with the line $g_1$.
Let it be a red line denoted by $r_1$; take its corresponding rotates $r_2$ and $r_0$.
\item \label{step:greenpoints} 
Take the point of intersection $G_0 \coloneqq p_0 \cap r_2$, and its rotates.
\item \label{step:yellowlines}
Take the circumcircle $\gamma_Y$ of the points $Y_0,Y_1,Y_2$, and invert it in the circle $\Omega$. 
Denote the inverse circle by $\gamma_y$.
Draw tangents to $\gamma_y$ from the point $P_1$, and choose the one subtending the greater angle with the line $g_1$.
Let it be a yellow line denoted by $y_1$; take its corresponding rotates $y_2$ and $y_0$.
\item \label{step:closingstep}
Take the point of intersection $M_0 \coloneqq    p_2 \cap r_2$, and its rotates.
\end{enumerate}
\medskip

The construction of QC(B) is thus complete. It is directly seen that this structure is \emph{movable}: indeed, it can be transformed into infinitely many 
projectively inequivalent versions while preserving incidences, by shifting the point $Y_2$ on the line $m_0$ (within the interval $R_1M_{01}$, as determined 
in step (\ref{step:yellowpoint})). Thus, the degree of freedom is 1; using the auxiliary elements in step (\ref{step:aux}) serves the very purpose of preventing 
larger degree of freedom from appearing.

By an easy check one sees that each point has its dual counterpart and vice versa. The duality is induced by reciprocity with respect to the circle $\Omega$ 
constructed in step (\ref{step:mid}); thus, duality $(P)\!\leftrightarrow\!(p)$ is defined there. The orbit $(r)$ is defined in step (\ref{step:redlines}) using the 
same reciprocity, thus duality $(\mathrm R)\!\leftrightarrow\!(r)$ also holds. Similarly, step (\ref{step:yellowlines}) defines duality $(Y)\!\leftrightarrow\!(y)$. 
By comparing the definition of orbit (G) in step (\ref{step:greenpoints}) with the definition and a property of orbit (g) in steps (\ref{step:greenlines}) 
and (\ref{step:purplepoints}), and using the already established dualities, one sees that duality $(G)\!\leftrightarrow\!(g)$ holds as well. Finally, duality 
$(M)\!\leftrightarrow\!(m)$ can be verified similarly by a comparison of steps (\ref{step:initialstep}) and (\ref{step:closingstep}). Observe that colouring 
the points and lines indicates their duality. 

From the type $((6_2)(9_4))$ of QC(B) one obtains that there are 48 incidences; due to symmetry, this means 16 incidence types, i.e.\ combinations  of 
colours of the form $(X,x')$. From Figure~\ref{fig:quasi} one can see that because of the presence of two trilaterals (with red vertices and with magenta 
vertices) this number reduces to 14. 13 types of incidences are defined directly in the construction; in the order of occurrence, these are: $(R,m)$, $(Y,m)$, 
$(R,p)$, $(Y,p)$, $(R,g)$, $(P,g)$, $(P,m)$, $(P,r)$, $(G,p)$, $(G,r)$, $(P,y)$, $(M,p)$, $(M,r)$. One observes that by the definition of the points and lines, 
and their duality established above, 12 of these incidences can be ordered into 6 dual pairs of the form $((X, x'),(X', x))$. Finally, the last, 13th incidence 
$(M,y)$ also holds, which is verified by duality since we have $(Y,m)$. 
\end{proof}


\subsection{Proof of Theorem~\ref{thm:MainThm}} \label{sect:MainThmProof}


\begin{proof}[Proof of Theorem~\ref{thm:MainThm}]
Here we continue the construction given above with additional steps. 

\begin{enumerate} [(1)]
\addtocounter{enumi}{11}
\item 
Define $B_0$ as the point of intersection $B_0\coloneqq g_0 \cap y_2$, and similarly its rotates by cyclic permutation of the indices.
\item \label{step:cyanlines} 
Define the cyan line $c_0\coloneqq   B_0B_2$, and similarly its rotates by cyclic permutation of the indices.
\item \label{step:bluelines} 
Take the circumcircle $\gamma_B$ of the blue points $B_0$, $B_1$, $B_2$, and invert it with respect to the circle $\Omega$. Denote the resulted circle 
by $\gamma_b$. Draw tangents to this circle from the point $Y_2$, and choose the one that is separated from the centre $O$ by the line $p_2$.
Denote it $b_1$, and take its rotates $b_2$ and $b_0$.
\item 
Define $C_0$ as the point of intersection $C_0\coloneqq b_0 \cap b_1$, and similarly its rotates by cyclic per\-mu\-tation of the indices.
\item 
Take the point of intersection of $c_1$ and $m_0$, and denote it by $Y'_2$. 

Observe that the incidence structure obtained in steps 1--16 above is movable. In fact, QC(B) is movable, as we have seen in the preceding subsection,
and this poperty has been preserved throughout the additional steps 12--16. In the following, we utilize this property. Indeed, consider the mutual position 
of the points $Y_2$ and $Y'_2$ within the interval $(R_1M_{01})$. 
\item 
Move $Y_2$ along the given interval. 
\end{enumerate}

Now, one observes that while shifting $Y_2$ along this interval, $Y'_2$ moves in the opposite direction. In particular, both 
cases of the order of the following four points may occur, such as $(R_1, Y'_2, Y_2, M_{01})$ and $(R_1, Y_2, Y'_2, M_{01})$ (this is checked in a 
dynamic geometry model). Hence, by continuity, an intermediate case must exist, where $Y_2$ and $Y'_2$ will coincide. This is precisely the case 
where our movable structure obtained above becomes equal to $\mathrm B(21_4)$.

This continuity argument is justified by the following observation. When building the incidence structure through steps (1--16), in each step where 
a new point or line is defined, its position is a continuous function of that of the already existing geometric constituents used in the definition. This 
is in fact a simple observation, provided we restrict ourselves to the domain $(R_1M_{01})$, where the first movable point $Y_2$ is located. 

Now we check the new incidences which occurred in this second part of our construction. Incidences $(B,g)$, $(B,y)$ and $(B,c)$ exist by definition.
Since orbit $(b)$ is defined by reciprocating $(B)$, the dual incidences $(G,b)$, $(Y,b)$ also hold. Again, we have $(C,b)$ by definition. It follows
the triangles $B_0B_1B_2$ and $C_0C_1C_2$ are dual to each other (by reciprocating in $\Omega$). This means, in particular, that orbits $(C)$
and $(c)$ are dual to each other. From the continuity argument we infer that incidence $(Y,c)$ holds. By duality, $(C,y)$ also holds. 

\begin{figure}[!h]
\begin{center}
\includegraphics[width=0.825\textwidth]{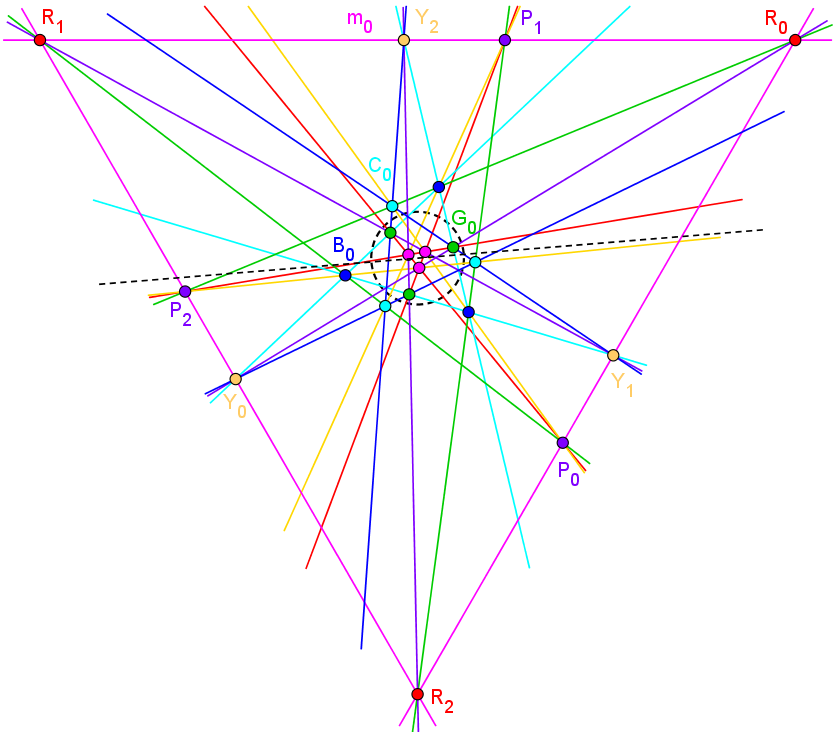}
\caption{The configuration $\mathrm B(21_4)$ with the circle of reciprocity indicated. 
The mirror line belonging to the self-reciprocation $\mathbf d$ given in Section~\ref{sect:aut} 
is also shown by a dashed line.}
\label{fig:selfrec}
\end{center}
\end{figure}

What remain to be verified are incidences $(C, g)$ and $(G,c)$; in fact, either of them is sufficient by duality. We choose $(C,g)$. 

Draw the auxiliary circle $\Gamma_a$ through the points $Y_0$, $O$, $R_1$. Due to the rotational symmetry, the angles $R_0Y_2O$ and $R_1Y_0O$ 
are equal. On the other hand, the angles $R_0Y_2O$ and $OY_2R_1$ are supplementary, thus so are the angles $OY_2R_1$ and $R_1Y_0O$. Hence 
the quadrangle $R_1Y_0OY_2$ is cyclic, which means that its vertex $Y_2$ also lies on the circumference of $\Gamma_a$. 

Since its angle at vertex $R_1$ is $60^{\circ}$, the angle at the opposite vertex $O$ is $120^{\circ}$. Observe 
that this latter angle is subtended by the arc $Y_2R_1Y_0$. In addition, the angle $Y_0C_1Y_2$ is subtended by the same arc, and (since it is the point 
of intersection of two blue lines) it is also $120^{\circ}$; hence its vertex $C_1$ lies on the circumference of $\Gamma_a$.

Consider now the angles $OR_1C_1$ and $OY_2C_1$. Since they are subtended by the same arc (determined by $C_1$ and $O$), while their vertex 
$R_1$ resp.\ $Y_2$ lie on the circumference of $\Gamma_a$, they are equal. Take the vertices of the triangle determined by the green lines, and denote
them by $U_0$, $U_1$ and $U_2$. Observe that the line $OU_1$ bisects the angle at $U_1$ of the triangle $U_0U_1U_2$, and so does the line $OR_1$
the angle at $R_1$ of the triangle $R_0R_1R_2$; hence both of these angles are $30^{\circ}$. It follows that the triangles $OR_1U_1$ and $OY_2C_1$ 
are similar to each other (these triangles are highlighted in Figure~\ref{fig:simproof} by red and yellow, respectively).

Recall that by a basic theorem of Euclidean plane geometry, any two similar triangles $A_1A_2A_3$ and $A_1'A_2'A_3'$ determine a unique similarity 
transformation which sends the vertex $A_i$ to vertex $A_i'$ for all $i\in\{1,2,3\}$. If the triangles are oriented alike, then this transformation is a dilative 
rotation~\cite{Cox1969}. In our case, we see that we have a dilative rotation $\varrho$ such that its fixed point is the common vertex $O$ of the red and 
yellow triangle considered above, and it acts on the two other vertices as follows: 
$$
\varrho: R_1 \mapsto Y_2,\; U_1 \mapsto C_1.
$$ 
Taking into account the threefold rotational symmetry, we see that $\varrho$ transforms the equilateral triangle $R_0R_1R_2$ into the triangle $Y_1Y_2Y_0$
inscribed in it (in the sense that the vertices of the latter are incident to the side lines of the former). Thus we conclude that $\varrho$ acts on the triangle
$U_0U_1U_2$ in the same way, which means that this triangle is transformed into the triangle whose side lines are the blue lines, vertices are the cyan 
points, moreover, it is inscribed in the triangle $U_0U_1U_2$. But the side lines of the triangle $U_0U_1U_2$ are the green lines, hence we see that the 
incidence  $(C,g)$ holds indeed.
\end{proof}

\begin{figure}[!h]
\begin{center}
\includegraphics[width=0.675\textwidth]{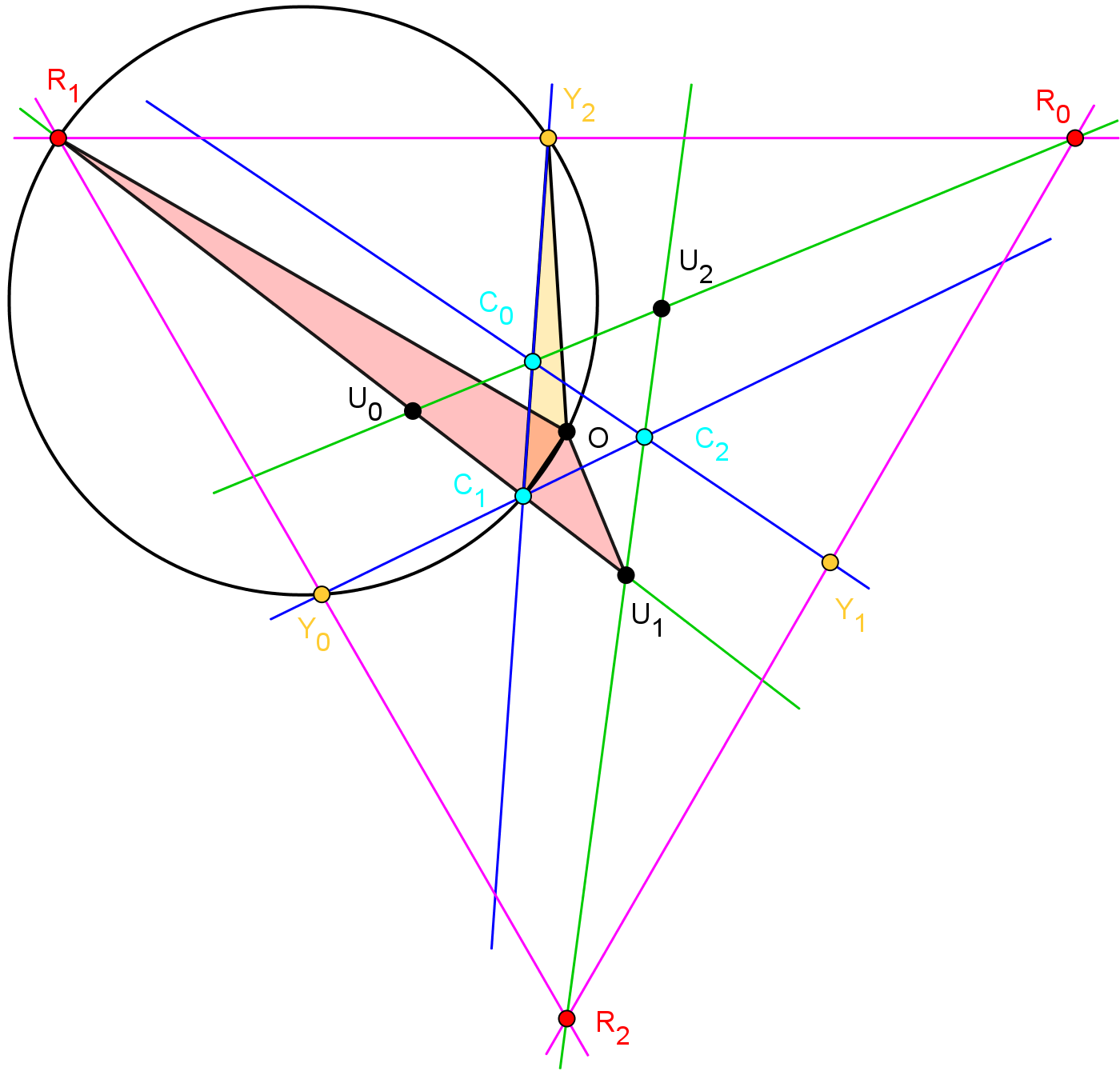}
\caption{Illustration for the verification of incidence $(C,g)$.}
\label{fig:simproof}
\end{center}
\end{figure}

Based on the construction above, here we provide the coordinates of the initial four points determining the configuration. Assume that the red points are located as follows:
\begin{equation*}
R_0 = (\sin(120^{\circ}), 0.5); \quad R_1 = (\sin(-120^{\circ}), 0.5); \quad R_2= (0, -1).
\end{equation*}
Then $Y_2$ can be given (up to 15 decimals) as follows:
\begin{equation*}
Y_2= (-0.031440363334572, 0.5).
\end{equation*}


\section{Automorphims and self-dualities of $\mathrm B(21_4)$} \label{sect:aut}

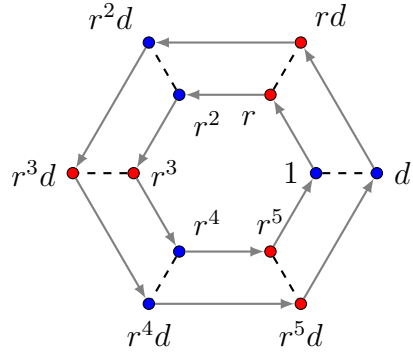
\begin{figure}[!h]
\begin{center}
\begin{tikzpicture}[scale = .8, vtx/.style={draw, circle, inner sep = 1.5 pt, font = \tiny}, lbl/.style={midway, inner sep = 1 pt, fill =white, font = \footnotesize}, lin/.style={draw, square, inner sep = 12 pt}]
\pgfmathsetmacro{\r}{1.5}
\pgfmathsetmacro{\rr}{2.5}
\node[fill=blue, vtx, label=left:{$1$}] (1) at (360*0/3: \r) {};
\node[fill=blue, vtx, label=-60:{$r^{2}$}] (r2) at (360*1/3:\r) {};
\node[fill=blue, vtx, label=60:{$r^{4}$}] (r4) at (360*2/3:\r) {};
\node[fill=blue, vtx, label=0:{$d$}] (d) at (360*0/3: \rr) {};
\node[fill=blue, vtx, label=90+60:{$r^{2}d$}] (r2d) at (360*1/3:\rr) {};
\node[fill=blue, vtx, label=below:{$r^{4}d$}] (r4d) at (360*2/3:\rr) {};

\node[fill=red, vtx, label=180+60:{$r$}] (r) at (360*0/3+60: \r) {};
\node[fill=red, vtx, label=right:{$r^3$}] (r3) at (360*1/3+60:\r) {};
\node[fill=red, vtx, label=90+1:{$r^{5}$}] (r5) at (360*2/3+60:\r) {};
\node[fill=red, vtx, label=60:{$rd$}] (rd) at (360*0/3+60: \rr) {};
\node[fill=red, vtx, label=left:{$r^{3}d$}] (r3d) at (360*1/3+60:\rr) {};
\node[fill=red, vtx, label=below:{$r^{5}d$}] (r5d) at (360*2/3+60:\rr) {};

\foreach \i/\j in {1/r,r/r2,r2/r3,r3/r4,r4/r5, r5/1, d/rd,rd/r2d,r2d/r3d,r3d/r4d,r4d/r5d, r5d/d}{
\draw[thick, gray, ->] (\i)--(\j);}

\foreach \i/\j in {1/d, r/rd,r2/r2d,r3/r3d,r4/r4d,r5/r5d}{\draw[dashed,  thick] (\i)--(\j);}

\end{tikzpicture}

\caption{
The Cayley graph of $\mathrm {Aut}( L(B))$ corresponding to the presentation (\ref{eq:presentation}); the generator $r$ is shown with gray arrows and the involutory generator $d$ with dashed edges.
The red vertices represent purely combinatorial automorphisms which cannot be realized geometrically (that is, using reciprocity and reflections or rotations).}
\label{fig:Cayley}
\end{center}
\end{figure}

We denote the automorphism group of the Levi graph of $\mathrm B(21_4)$ by $\mathrm{Aut} \mathrm{L(B)}$. 
It is isomorphic to the dihedral group $\mathrm D_6$ of order 12 (see item 1 in Table~\ref{tab:42cases}).
As usual, it is generated by two generators $\mathbf r$, $\mathbf d$ together with the following defining relations: 
\begin{equation} \label{eq:presentation}
\mathbf r^6 = \mathbf d^2 = (\mathbf {rd})^2 =\mathbf 1.
\end{equation}
The corresponding Cayley graph of this group is depicted in Figure~\ref{fig:Cayley}.

In our case, the generators take the following form:

\begin{equation} \label{eq:order6}
\begin{split}
{\bf r}=&(R_0,M_0,R_1,M_1,R_2,M_2)(Y_0,G_2,Y_1,G_0,Y_2,G_1)\\
            &(B_0,C_2,B_1,C_0,B_2,C_1)(P_0,P_2,P_1)\\
            &(r_0,m_1,r_1,m_2,r_2,m_0)(y_0,g_2,y_1,g_0,y_2,g_1)\\
            &(b_0,c_2,b_1,c_0,b_2,c_1)(p_0,p_2,p_1),\\
\end{split}
\end{equation}
\begin{equation} \label{eq:selfrec}
\begin{split}
\mathbf d=&(R_0,r_1)(R_1,r_0)(R_2,r_2)(Y_0,y_1)(Y_1,y_0)(Y_2,y_2)\\
              &(P_0,p_1)(P_1,p_0)(P_2,p_2)(G_0,g_1)(G_1,g_0)(G_2,g_2)\\
              &(B_0,b_1)(B_1,b_0)(B_2,b_2)(C_0,c_1)(C_1,c_0)(C_2,c_2)\\
              &(M_0,m_1)(M_1,m_0)(M_2,m_2).
\end{split}
\end{equation}

Using Table~\ref{tab:inci}, one can directly check that ${\bf r}$ given by equation (\ref{eq:order6}) is indeed an automorphism.

On the other hand, taken to the second power
\begin{equation}
\begin{split}
{\bf r}^2=&(R_0,R_1,R_2)(r_0,r_1,r_2)(Y_0,Y_1,Y_2)(y_0,y_1,y_2)(P_0,P_1,P_2)(p_0,p_1,p_2)\\
                &(G_0,G_1,G_2)(g_0,g_1,g_2)(B_0, B_1, B_2)(b_0,b_1, b_2)(C_0,C_1,C_2)(c_0,c_1,c_2)\\
                &(M_0,M_1,M_2)(m_0,m_1,m_2),
\end{split}
\end{equation}
one obtains a (geometric) rotation of order 3, which generates the symmetry group 
$\mathrm{Sym}\,\mathrm B(21_4)$ $=\langle \mathbf r^2 \rangle \cong \mathbb Z_3$ 
of $\mathrm B(21_4)$. We note that this verifies as well that $\mathrm B(21_4)$ is a polycyclic configuration. 

The automorphism ${\bf r}$ generates the full automorphism group $\mathrm{Aut}\,\mathrm B(21_4) \cong \mathbb Z_6$ of our 
configuration. The other 6 elements of $\mathrm{Aut} \mathrm{L(B)}$ are all (involutory) self-dualities of $\mathrm B(21_4) $. 
Three of them are purely combinatorial self-dualities (see Figure \ref{fig:Cayley}). On the other hand, $\mathbf d$ given by equality 
(\ref{eq:selfrec}), and its two conjugates $\mathbf r^2\mathbf d$ and $\mathbf r^4\mathbf d$ are realized as geometric transformations. 
Namely, they are reflexible self-reciprocations of $\mathrm B(21_4) $. The circle of reciprocity is shown in Figure~\ref{fig:selfrec}. The 
mirror line belonging to $\mathbf d$ is also shown in the same figure; the mirror lines belonging to the two other self-reciprocations are 
rotates of this line by angles $\pm 120^{\circ}$.

As it can directly be seen from Figure \ref{fig:Cayley}, one of the purely combinatorial self-dualities is given as ${\bf r}^5{\bf d}={\bf dr}$.
Using equalities (\ref{eq:order6}) and (\ref{eq:selfrec}), for this product one obtains the following expression: 
\begin{equation} \label{eq:combinatorialdual}
\begin{split}
\mathbf{dr}=&(R_0,m_2)(R_1,m_1)(R_2,m_0)(Y_0,g_0)(Y_1,g_2)(Y_2,g_1)\\
              &(P_0,p_0)(P_1,p_2)(P_2,p_1)(G_0,y_0)(G_1,y_2)(G_2,y_1)\\
              &(B_0,c_0)(B_1,c_2)(B_2,c_1)(C_0,b_0)(C_1,b_2)(C_2,b_1)\\
              &(M_0,r_1)(M_1,r_0)(M_2,r_2).
\end{split}
\end{equation}

Observe that the symmetry properties of the reduced Levi graph of $\mathrm B(21_4) $ show both types of self-duality: indeed, its half-turn symmetry 
shows self-reciprocity, while the mirror symmetry with respect to a horizontal axis shows precisely the type of combinatorial duality given by the expression 
above (see the colours of the nodes representing the point orbits and line orbits of the configuration).

We conclude this subsection with some questions related to the rank of self-duality. Let $\mathcal C$ be a self-dual configuration, and let $\delta$ denote a
self-duality map of $\mathcal C$. The \emph{rank} $r(\delta)$ of $\delta$ is defined as its order, i.e.\ the smallest positive integer $n$ such that $\delta^n$ is 
the identity. A configuration may have more than one self-duality maps (in our case we have altogether 6, as we have seen above), thus the following definition 
makes sense. The \emph{rank} $r(\mathcal C)$ of a self-dual configuration $\mathcal C$ is the minimum value of $r(\delta)$ over all self-duality maps $\delta$
of $\mathcal C$. (We note that this notion was introduced by Gr\"unbaum and Shephard~\cite{GruShe1988} in case of geometric objects for which self-duality 
can be defined, in particular, for polyhedra and configurations; for further details related to configurations, see~\cite{Gru2009b} and the references therein).

As we have seen above, all the self-dualities of $\mathrm B(21_4)$ are involutory, thus its rank $r(\mathrm B(21_4))= 2$. But as we have also seen, 
$\mathrm B(21_4)$  is reflexibly self-reciprocal; hence the question arises that, in general, does the latter property imply the former? It is appropriate 
here to cite a conjecture by Gr\"unbaum~\cite[Conjecture 5.8.1]{Gru2009b}.

\begin{conjecture}
Every self-dual geometric configuration of rank 2 has a realization such that its polar (in a suitable circle) is congruent to the original configuration.
\end{conjecture}


\section{Analytic proof of the existence of the configuration $B(21_{4})$} \label{sect:analytic}


We use the reduced Levi graph shown in Figure \ref{fig:RLG(B)} as a ``recipe'' to construct the configuration analytically with homogeneous coordinates, 
using \emph{Mathematica}. We are interested in constructing a strong realization of the configuration, in which all the points and lines are distinct.

There are a number of ways to walk through the reduced Levi graph. Here, we present one which uses only meets and joins so that in the final steps, 
we are solving a system of two polynomials in two unknowns. We begin by fixing the red points and magenta lines. (For convenience, we use 
a rotated starting position for the red points from that used in the previous section.) We make the following assignments, following the labels on the 
reduced Levi graph shown in Figure \ref{fig:RLG-general} with the assignments from Figure \ref{fig:RLG(B)}. 
\[ \{a, c, d, e, f, g, q, a', c', d', e', f', g', q', t\}  =  \{1, 2, 1, 1,
   1, 2, 1, 1, 2, 1, 1, 1, 2, 1, 0\}.\]
 We place a yellow point $Y_{i}$ arbitrarily (using parameter $x$) on a magenta line 
$m_{i-a}$, and then use the red points $R_{i}$ and yellow points $Y_{i}$ to define the purple lines $p_{i}$. On the purple line $p_{i}$, we place a 
green point $G_{i}$ arbitrarily, using parameter $z$. The rest of the points and lines are determined, as follows, in order (all indices taken mod 3):
\begin{align}\label{eqn:pathThruRLG}
\begin{split}
R_{i} &= (2\cos(2 \pi i/3), 2 \sin(2 \pi i/3), 1)\\
 m_{i} &= R_{i} \times R_{i+a}\\
 Y_{i}&= (1-x) m_{i-c} + x m_{i-c+a}\\
 p_{i}&= Y_{i} \times R_{i}\\
G_{i}&= (1-z) R_{i} + z Y_{i}\\
 b_{i} &= Y_{i+e'} \times G_{i}\\
 c_{i} &= Y_{i} \times G_{i+f'}\\
C_{i} &= b_{i-g} \times b_{i-t}\\
 B_{i} &= c_{i-g'} \times c_{i}\\
g_{i} &= C_{i+f} \times B_{i}\\
 y_{i} &= C_{i}\times B_{i+e}\\
 M_{i}&= y_{i-c'} \times p_{i-q'}\\
 r_{i}&= G_{i+d} \times M_{i}
 \end{split}
\end{align}

As we made these assignments, we eliminated any common numeric or polynomial factors from the homogeneous coordinates.
After these simplifications, we computed two determinants: the reduced Levi graph says that we need $R_{d'}, C_{f}, B_{0}$ 
collinear (on $g_{0}$) and $M_{0}, M_{a'}, G_{d}$ collinear (on $r_{0}$). Define
\begin{align*} 
{\tt det1} &= \det( R_{d'} \ C_{f} \ B_{0} ) \\
&= -6 \sqrt{3} (x-1) z \left(2 x^2 z+x^2-2 x z-x+z\right) \left(3 x^3 z-x^2 z^2-4 x^2 z-x^2+x z^2+3 xz-z^2\right) \end{align*}
and
\begin{align*}
{\tt det5} &= 
\det( M_{0} \ M_{a'} \ G_{d} ) \\
&= 6 \sqrt{3} \left(x^2 z-x z-x+z\right) \left(2 x^2
   z+x^2-2 x z-x+z\right) \left(3 x^6 z^3+6 x^6 z^2-7
   x^5 \right.\\& \phantom{=} \qquad \left. z^3-27 x^5 z^2-2 x^5+9 x^4 z^3+53 x^4 z^2-9 x^4
   z+7 x^4-6 x^3 z^3-62 x^3 z^2+
   \right.\\& \phantom{=} \qquad \left.
   25 x^3 z-11 x^3+2 x^2
   z^3+ 44 x^2 z^2-28 x^2 z+10 x^2-18 x z^2+15 x z-5 x+3
   z^2-3 z+1\right)
   \end{align*}
   
These two determinants have a common factor 
\[{\tt common} = 6 \sqrt{3} \left(2 x^2 z+x^2-2 x z-x+z\right); \]
solving ${\tt common} = 0$ for $z$ in terms of $x$ yields 
\begin{equation}\label{eqn:common} z = \frac{x-x^2}{2 x^2-2 x+1}.\end{equation}
Performing this substitution into both $M_{-q'}$ and $Y_{0}$, for example, which both lie on the purple line $p_{0}$, shows that the two points coincide, 
which is forbidden in a strong realization of the reduced Levi graph.

Reducing the two determinants into ${\tt det1'}$ and ${\tt det5'}$ respectively by eliminating the common factor and solving the associated system 
\begin{equation}\label{eqn:det}
\{ {\tt det1'}=0,  {\tt det5'}=0\}
\end{equation} over $\mathbb{R}$ results in a collection of solutions. Some of the solutions are degenerate,
\[\left\{\{x\to 0,z\to 0\right\},\left\{x\to \frac{1}{2},z\to 0\right\},\left\{x\to \frac{1}{2},z\to \frac{2}{3}\right\},\{x\to 1,z\to 0\},\{x\to 1,z\to 1\},\]because 
these all immediately lead to coinciding points in the construction.

However, there are exactly two non-degenerate solutions over $\mathbb{R}$, which \emph{Mathematica} can express exactly as roots of certain polynomials with integer 
coefficients: let 
\begin{align*}
\alpha(s) &= 9 s^{12}-45 s^{11}+108 s^{10}-114 s^9-57 s^8+\\
& \qquad 390
   s^7-668 s^6+684 s^5-468 s^4+217 s^3-66 s^2+12 s-1\\
\beta(s) &=  s^{12}+9 s^{11}+15 s^{10}-36 s^9-33 s^8+\\
& \qquad 129 s^7-193
   s^6+216 s^5-162 s^4+76 s^3-24 s^2+6 s-1  
   \end{align*}
   
 Both of these polynomials have exactly two real roots. The two solutions $(x, z)$ to \eqref{eqn:det} are

\begin{align*} x &= \text{Root}(\alpha(s), 1)  \approx -1.66271\\
   z &= \text{Root}(\beta(s), 1) \approx -5.40326
   \end{align*}
\noindent   
and 
 \begin{align*} x &= \text{Root}(\alpha(s), 2)  \approx 0.518152\\
   z &= \text{Root}(\beta(s), 2) \approx 0.611257
   \end{align*}  
 \noindent  
corresponding to whether the magenta points are outside or inside the circumcircle of the red points. It is straightforward to check via computer algebra that 
these solutions do not satisfy \eqref{eqn:common} and that the point sets  corresponding to these solutions are all distinct.

Finally, we need to show that the purple points lie on the intersections of four lines. To do so, we again compute two determinants corresponding to magenta, 
green and yellow lines concurrent, and red, green, and yellow lines concurrent, which turn out to be even higher-degree polynomials in $x$ and $z$:
\begin{align*}
{\tt det3} &= \det(m_{-q} \ g_{0} \ y_{0} )\\
&= -36 \sqrt{3} (z-1) \left(3 x^7 z^4+3 x^7 z^3+3 x^7 
   z^2-6 x^6 z^4-15 x^6 z^3-15 x^6 z^2   \right.\\ & \qquad \qquad \left.+5 x^5 z^4+31 x^5
   z^3+27 x^5 z^2-2 x^5 z-x^5
+2 x^4 z^4-40 x^4 z^3   \right.\\ & \qquad \qquad \left. -25
   x^4 z^2+6 x^4 z+3 x^4-8 x^3 z^4+34 x^3 z^3+13 x^3
   z^2-8 x^3 z-3 x^3+8 x^2 z^4\right.\\ & \qquad \qquad \left.-20 x^2 z^3-3 x^2 z^2+
   6x^2 z+x^2-4 x z^4+8 x z^3-x z^2-2 x z+z^4-2
   z^3+z^2\right)
\end{align*}
and
\begin{align*}
{\tt det4} &= \det(r_{0} \ g_{0} \ y_{0} )\\
&= -108 (z-1)^3 \left(3 x^{12} z^5+3 x^{12} z^4+3 x^{12}
   z^3-9 x^{11} z^5+9 x^{11} z^4\right.\\ & \qquad \left.+   9 x^{11} z^3+9 x^{11}
   z^2
+6 x^{10} z^5-108 x^{10} z^4-123 x^{10} z^3-63
   x^{10} z^2-3 x^{10} z\right.\\ & \qquad \left.+35 x^9 z^5+377 x^9 z^4+422 x^9
   z^3+182 x^9 z^2+11 x^9 z-x^9-129 x^8 z^5\right.\\ & \qquad \left.-782 x^8
   z^4-821 x^8 z^3-306 x^8 z^2-17 x^8 z+4 x^8+246 x^7
   z^5\right.\\ & \qquad \left.
+1112 x^7 z^4+1059 x^7 z^3+339 x^7 z^2+14 x^7 z-6
   x^7-318 x^6 z^5-1145 x^6 z^4\right.\\ & \qquad \left.-953 x^6 z^3-255 x^6
   z^2-5 x^6 z+4 x^6+300 x^5 z^5+866 x^5 z^4+602 x^5
   z^3+126 x^5 z^2-x^5 z\right.\\ & \qquad \left.
-x^5-210 x^4 z^5-476 x^4
   z^4-258 x^4 z^3-36 x^4 z^2+x^4 z\right.\\ & \qquad \left.+108 x^3 z^5+182 x^3
   z^4+68 x^3 z^3+4 x^3 z^2-39 x^2 z^5-44 x^2 z^4-8 x^2
   z^3+9 x z^5+5 x z^4-z^5\right)
   \end{align*}
   
Evaluating each of these determinants using the exact solutions found above (and the power of \emph{Mathematica's} symbolic algebra computations) shows that 
both these determinants evaluate to exactly 0, showing that the four lines $r_{0}, g_{0}, y_{0}$, and $m_{-q}$ are concurrent, at a point we label 
$P_{0}$ (and by symmetry, the other two quadruples of lines are concurrent at $P_{1}$, $P_{2}$).


\section{Polycyclic realizations of  the $(21_{4})$ Gr\"unbaum--Rigby configuration } \label{sect:gen}


The construction and investigation of $\mathrm{B}(21_4)$ can be generalized in several ways; here we outline some possibilites.

\subsection{Changing the voltage assignments in RLG(B)} \label{sect:voltages}

We systematically explored all voltage assignments that gave rise to a Levi graph of some combinatorial $(21_4)$ configuration over the  
quotient graph RLG(B) with generic parameters shown in  Figure~\ref{fig:RLG-general}. Two of the authors independently wrote 
computer programs (TP in \emph{Sage}; LWB in \emph{Mathematica}) that checked all possible values of the parameters \[\{a, c, d, e, f, g, q, a', c', d', e', f', g', q', t\}\]
over $\mathbb{Z}_{3}$, and eliminated parameter lists that produced Levi graphs 
with girth less than 6 and isomorphic graphs. We found 17 such Levi graphs. Their parameters are shown in Table \ref{tab:42cases}.  Line 1 in this table corresponds to RLG(B), 
while line 3 corresponds to 
the Gr\"unbaum-Rigby configuration with threefold (rather than seven-fold) rotational 
symmetry; that is, line 3 provides parameter values for the reduced Levi graph shown in Figure \ref{fig:RLG-general} for which the corresponding Levi graph is isomorphic to the 
Levi graph of the Gr\"{u}nbaum-Rigby configuration.

However, there is no \emph{strong realization} of the reduced Levi graph with those parameters: that is, every realization of the Gr\"{u}nbaum-Rigby configuration over RLG(B) 
has symmetry classes of points which coincide with each other, which we demonstrate by the following proposition.

\begin{proposition}
The Gr\"{u}nbaum-Rigby configuration admits only one strong realization as a polycyclic geometric configuration. 
\end{proposition}

The proof uses a computer but could, in principle, be determined by hand.

\begin{proof}
As we indicated above, there are only two non-isomorphic polycyclic reduced Levi graphs for the Gr\"unbaum-Rigby configuration. This gives rise to two polycyclic combinatorial realizations,
one with 7-fold symmetry and the other one with 3-fold symmetry.  The configuration with 7-fold symmetry is polycyclically geometrically realizable in essentially only one way
(see \cite{BerBer2014, Gru2009b}). A sketch of the argument is that any polycyclic realization of the reduced Levi graph given in \ref{fig:GR} is a celestial configuration with symbol 
$m\#(s_{1}, t_{1}; s_{2}, t_{2}; s_{3}, t_{3})$, which must satisfy the cosine condition and three other conditions, described in \cite[Theorem 3.7.1]{Gru2009b}. However, the only 
solutions to the cosine condition for $m = 7$ are cyclic permutations of $m\#(2,1;3,2;1,3)$ and $m\#(3,1;2,3,1,2)$, which produce congruent geometric configurations.

For 3-fold rotational symmetry, \emph{Mathematica}-based symbolic calculations show that the configuration is not geometrically realizable.
To see this, use the same pathway and point and line coordinate assignments through the reduced Levi graph described in equation \eqref{eqn:pathThruRLG}, only this time using the parameter assignments 
\[ \{a, c, d, e, f, g, q, a', c', d', e', f', g', q', t\} = \{1, 2, 0, 1, 1, 2, 2, 1, 2, 0, 1, 1, 2, 2, 0\}\]
given in Table \ref{tab:42cases} line 3.
In this case, similarly to the previous construction, we define
\begin{align*} 
{\tt det1} &= \det( R_{d'} \ C_{f} \ B_{0} ) \\
&=  \left(2 x^2 z+x^2-2 x z-x+z\right)
   \left(3 x^4 z^2-6 x^3 z^2-x^2 z^3+6 x^2 z^2\right.\\
& \phantom{=}\left. \qquad   -x^2+x z^3-3 x z^2+x-z^3+2 z^2-z\right)
 \end{align*}
 
and
\begin{align*}
{\tt det5} &= \det( M_{0} \ M_{a'} \ G_{d} ) \\
&= \left(12 x^{10} z^5+3 x^{10} z^4+12x^{10} z^3-60 x^9 z^5\right.\\
&\phantom{=}\left. \qquad -15 x^9 z^4-60 x^9z^3+158 x^8 z^5+35 x^8 z^4\right.\\
&\phantom{=}\left. \qquad +128 x^8 z^3-16 x^8z^2-4 x^8 z-4 x^8-272 x^7 z^5\right.\\
&\phantom{=}\left. \qquad -50 x^7 z^4-152x^7 z^3+64 x^7 z^2+16 x^7 z+16 x^7\right.\\
&\phantom{=}\left. \qquad +334 x^6z^5+44 x^6 z^4+96 x^6 z^3-112 x^6 z^2-32 x^6z\right.\\
&\phantom{=}\left. \qquad -24 x^6-302 x^5 z^5-20 x^5 z^4-8 x^5 z^3+112x^5 z^2\right.\\
&\phantom{=}\left. \qquad +40 x^5 z+16 x^5+203 x^4 z^5-4 x^4z^4-40 x^4 z^3\right.\\
&\phantom{=}\left. \qquad -72 x^4 z^2-28 x^4 z-4 x^4-100x^3 z^5+13 x^3 z^4\right.\\
&\phantom{=}\left. \qquad +36 x^3 z^3+32 x^3 z^2+8 x^3z+35 x^2 z^5-10 x^2 z^4\right.\\
&\phantom{=}\left. \qquad -16 x^2 z^3-8 x^2 z^2-8x z^5+4 x z^4+4 x z^3+z^5-z^4\right)
   \end{align*}

Solving the system $\{ {\tt det1} = 0, {\tt det5 = 0}\}$ leads to the five real solutions 
\[ \{x = 0, z = 0\},\{x = 0, z =1\},  \{x = 1, z = 0\}, \{x = 1, z = 1\}, \left\{x = \frac{1}{2}, z = \frac{2}{3}\right\},\]
but all of these solutions lead to coinciding points and lines, and thus only a degenerate geometric realization of the configuration.

\end{proof}

\begin{theorem}
The configuration $B(21_{4})$ is the only $(21_{4})$ configuration with a nondegenerate $\mathbb{Z}_{3}$ geometric polycyclic realization. 
\end{theorem}

\begin{figure}[htbp]
\begin{center}

\includegraphics[width = .4\linewidth]{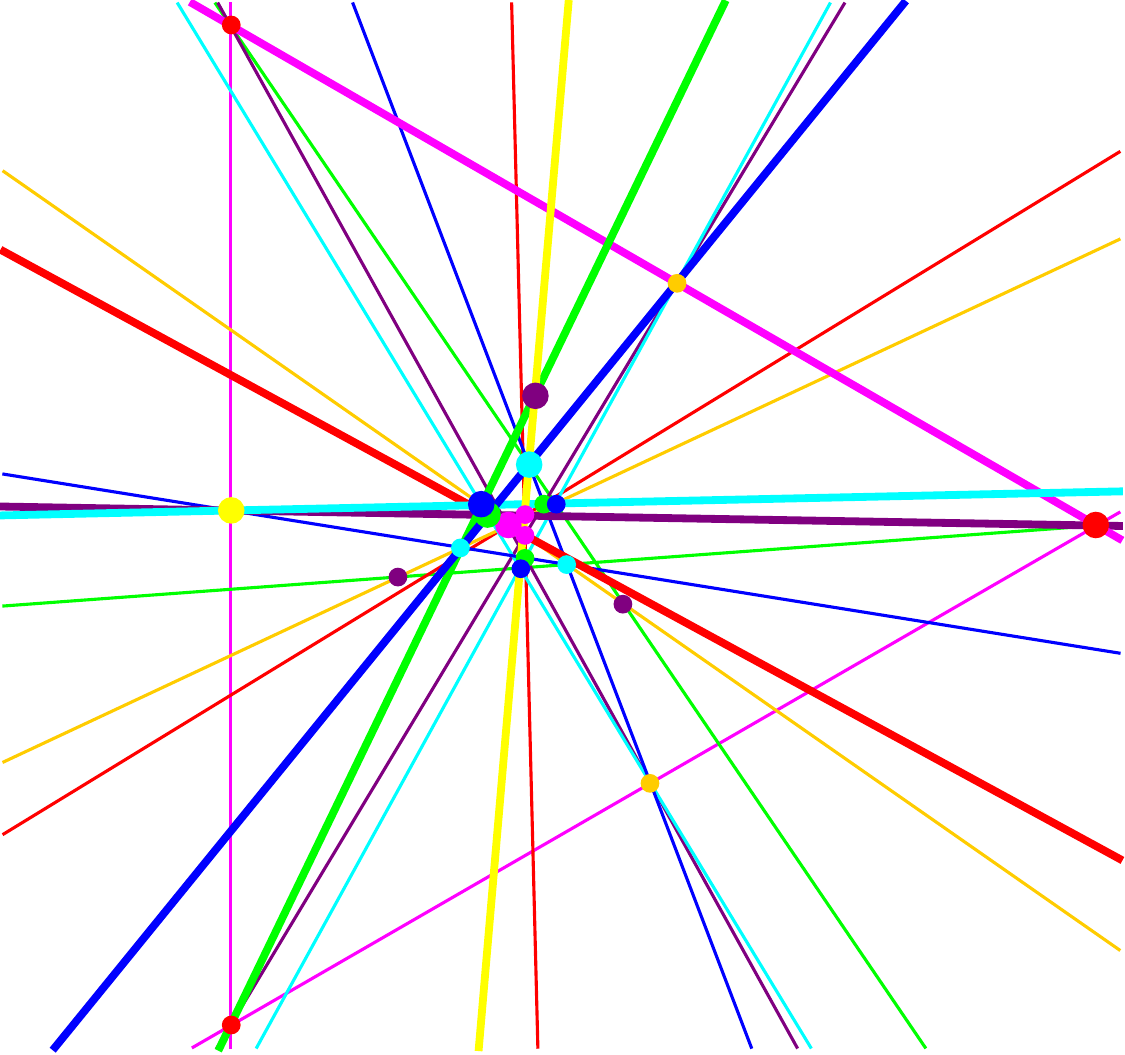} \hspace{1cm} 
\includegraphics[width = .4\linewidth]{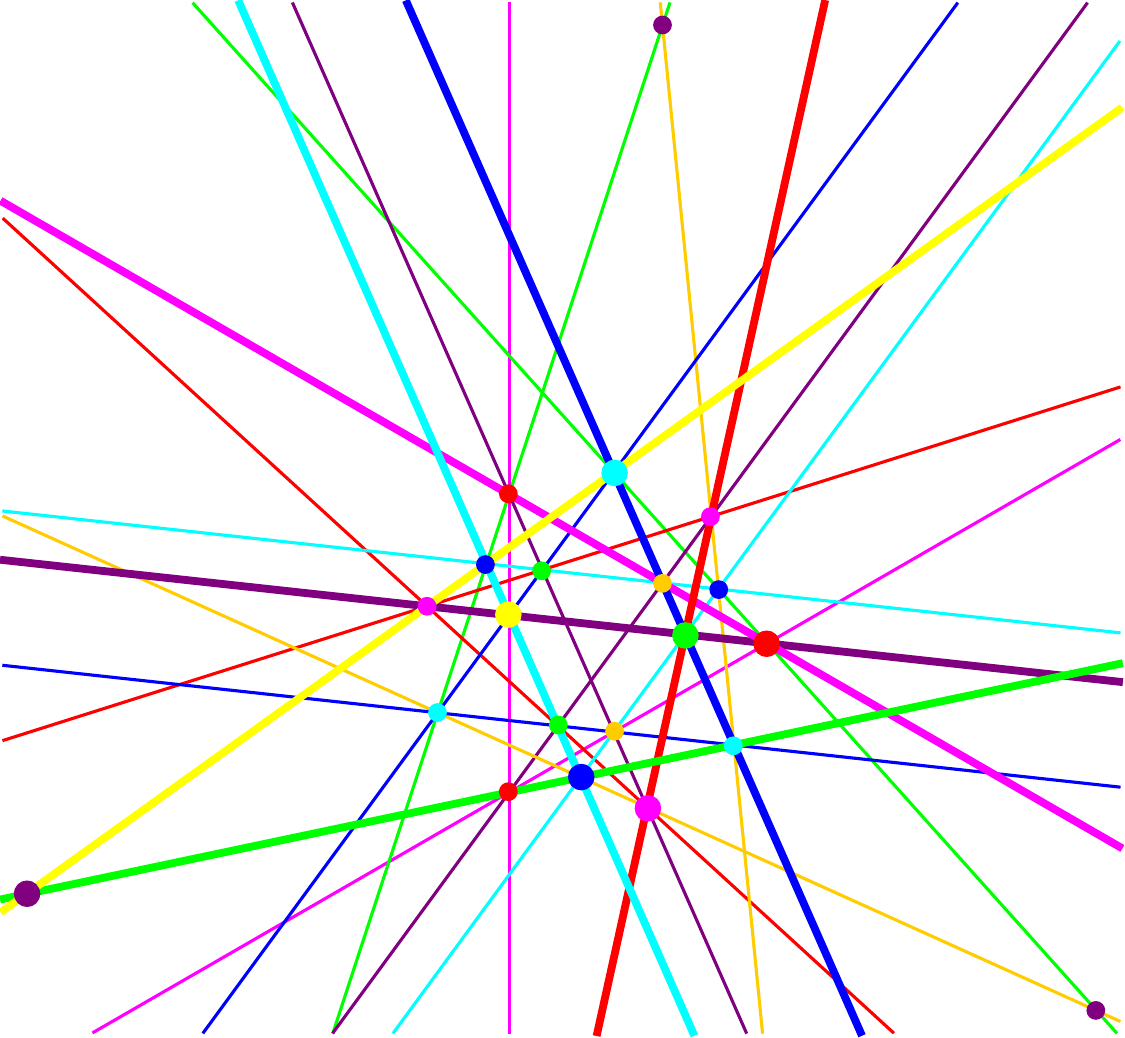}

\caption{Configurations \#9 (left) and \#13 (right) have at least one nondegenerate realization that satisfies Equation \eqref{eq:system}, which means that points  $R_{d'}, C_{f}, B_{0}$ 
are collinear (on the thick green line) and points $M_{0}, M_{a'}, G_{d}$ are collinear (on the thick red line). However, because the determinants in equation \eqref{eq:otherDets} are 
not identically equal to 0, the purple points do not lie on the common intersection of four lines. In these drawings, the purple points $P_{i}$ are defined as the intersection of lines $y_{i}$ 
and $g_{i}$, and the 0th element of each class is shown large (for points) or thick (for lines).}
\label{fig:badParamPix}
\end{center}
\end{figure}

\begin{proof}
As in the previous proposition, we analyzed all 17 parameter lists given in Table \ref{tab:42cases}, using the same point and line coordinate assignments through the reduced Levi graph 
described in equation \eqref{eqn:pathThruRLG} for each set of parameters. In each case, we assigned ${\tt det1} = \det( R_{d'} \ C_{f} \ B_{0} )$ and ${\tt det5} = \det( M_{0} \ M_{a'} \ G_{d} )$ 
and found all solutions to the system 
\begin{equation}\label{eq:system}\{ {\tt det1} = 0, {\tt det5 = 0}\}.\end{equation}

The $B$ configuration \#1 and the GR configuration \#3 have been analyzed above. Of the remaining configurations, \#4, \#5, \#6, \#8, \#10, \#11, \#12, \#15, \#16 only have degenerate 
solutions to \eqref{eq:system} (that is, all solutions lead to coinciding sets of points and lines).

For the remaining configurations \#2, 7, 9, 13, 14, 17, as in the analysis of $B$, we then evaluated the two determinants 
\begin{equation}
{\tt det3} = \det(m_{-q} \ g_{0} \ y_{0} ) \text{ and } {\tt det4} = \det(r_{0} \ g_{0} \ y_{0} ) \label{eq:otherDets}
\end{equation}
 at the nondegenerate solutions to equation \eqref{eq:system}. These two determinants must both evaluate to exactly 0 in order for the four lines $m_{i-q}$, $g_{i}$, $y_{i}$, $r_{i}$ 
 to pass through each purple point $P_{i}$.

The determinants in equation \eqref{eq:otherDets} for configurations \#2, 7, 14, 17 evaluated to numbers that were very far from 0 (on the order of $10^{7}$). In contrast, the values 
of the determinants for \#9, \#13 were numerically both between 0 and 1; however, computing the values of the determinants exactly showed (eventually) that they were not identically 
equal to 0 and thus, there is no nondegenerate geometric polycylic realization of either configuration. Pictures of both of these configurations are shown in Figure \ref{fig:badParamPix}.
\end{proof}

\begin{corollary} There are exactly two geometric polycyclic $(21_{4})$ configurations.\end{corollary}

\begin{proof} The configuration $B(21_{4})$ can be polycylically geometrically realized  with $\mathbb{Z}_{3}$ symmetry, and the configuration $GR(21_{4})$ can be polycylically 
geometrically realized  with $\mathbb{Z}_{7}$ symmetry. Since the only two possible reduced Levi graphs have each been analyzed and there are no other parameter values that 
lead to nondegenerate realizations, it follows that these are the only two geometric polycyclic configurations.
\end{proof}

\begin{table} [!t]
\begin{center}
\caption{Parameters of Levi graphs of all 17 combinatorial configurations derived from RLG(B), with parameters corresponding to those in~\ref{fig:RLG-general}. 
Lines 1 and 3 (highlighted in gray) correspond to $B$ and $GR$ respectively. The column |Aut| gives the number of automorphisms of the Levi graph. 
The column ``self-dual?'' indicates whether the configuration is combinatorially self-dual, and the column ``NDSols'' indicates whether there are any non-degenerate 
solutions to the system $\{ {\tt det1} = 0, {\tt det5 = 0}\}$ (see text); the annotation `y' says that there are non-degenerate solutions that do not lead to a full geometric realization, while `Y' says that there are non-degenerate solutions that {\bf do} lead to a full geometric realization.}
\begin{small}
\begin{tabular}{crcccc}
\hline
 item & $\{a,c,d,e,f,g,q,a',c',d',e',f',g',q',t\}$ & |Aut| & name & self-dual? & NDsols?\\
\hline
\rowcolor{gray!20}
1&$\{1, 2, 0, 1, 1, 1, 2, 1, 2, 0, 1, 1, 1, 2, 0\}$ & 12 & B & y & Y\\
2& $\{1, 2, 0, 1, 1, 1, 2, 1, 2, 0, 1, 1, 2, 2, 0\}$ & 6 & &y & y\\
\rowcolor{gray!20}
3&$\{1, 2, 0, 1, 1, 2, 2, 1, 2, 0, 1, 1, 2, 2, 0\}$& 672& GR & y & n\\
4&$\{1, 2, 0, 1, 1, 1, 2, 1, 2, 0, 2, 2, 1, 2, 0\}$& 12 & & y &n\\
5&$\{1, 2, 0, 1, 1, 1, 2, 1, 2, 0, 2, 2, 2, 2, 0\}$& 12& &y & n\\
6&$\{1, 2, 0, 1, 1, 1, 0, 1, 2, 2, 1, 1, 1, 2, 0\}$& 6 && y & n\\
7&$\{1, 2, 0, 1, 1, 1, 0, 1, 2, 2, 1, 1, 2, 2, 0\}$& 6&&y& y\\
8&$\{1, 2, 0, 1, 1, 1, 0, 1, 2, 2, 2, 2, 1, 2, 0\}$& 12&&y& n\\
9&$\{1, 2, 0, 1, 1, 1, 0, 1, 2, 2, 2, 2, 2, 2, 0\}$& 3 &&n& y\\
10&$\{1, 2, 0, 1, 1, 2, 0, 1, 2, 2, 2, 2, 2, 2, 0\}$& 6&&y&n\\
11&$\{1, 2, 0, 2, 2, 1, 0, 1, 2, 2, 2, 2, 2, 2, 0\}$& 6&&y&n\\
12&$\{1, 2, 0, 2, 2, 2, 0, 1, 2, 2, 2, 2, 1, 2, 0\}$& 6&&y&n\\
13&$\{1, 2, 0, 2, 2, 2, 0, 1, 2, 2, 2, 2, 2, 2, 0\}$& 12&&y&y\\
14&$\{1, 2, 2, 1, 1, 1, 0, 1, 2, 2, 2, 2, 2, 0, 0\}$& 6&&y&y\\
15&$\{1, 2, 2, 1, 1, 2, 0, 1, 2, 2, 2, 2, 2, 0, 0\}$& 24&&y&n\\
16&$\{1, 2, 2, 2, 2, 1, 0, 1, 2, 2, 2, 2, 2, 0, 0\}$& 12&&y&n\\
17&$\{1, 2, 1, 1, 1, 1, 0, 2, 1, 2, 2, 2, 2, 0, 0\}$& 6&&y&y\\
 \hline
\end{tabular}
\label{tab:42cases}
\end{small}
\end{center}
\end{table}

It is interesting that among the 17 non-isomorphic Levi graphs, 16 give rise to self-dual combinatorial configurations. Only one Levi graph, defined by parameters in Table \ref{tab:42cases} line \#9, gives 
rise to a pair of dual configurations, bringing the total of non-isomorphic configurations to 18. A Levi graph admits a self-dual configuration if and only if it has an automorphism 
that interchanges the sets of bipartition. One would expect that in such a situation, the dual pair of configurations would give rise to two sets of parameters. However, this is 
not the case here. Namely, the dual pair of configurations have isomorphic underlying Levi graphs (with vertex colors reversed). On the other hand, there is no color preserving 
isomorphism that would map one (vertex-colored) Levi graph onto the other one. The opposite is true in all other 16 cases. 


\section{Conclusions and open questions}


Since two non-isomorphic geometric $(21_4)$ configurations exist, and these are the only (strongly) realizable polycyclic geometric $(21_{4})$ configurations, the natural question is: 
\emph{Are there more of them?} An over-ambitious 
project involves a solution to the following formal problem.

\begin{problem}
Determine all geometric $(21_4)$ configurations.
\end{problem}

The complete solution to this problem seems to be out of reach with our current knowledge about configurations.
The brute-force approach does not seem feasible. Namely, no one knows how many combinatorial $(21_4)$ configurations
exist. It is not known how many connected bipartite graphs of girth at least 6 on 42 vertices exist. The number must be 
large, since is known that there exist almost two billion distinct quartic graphs of girth at least 6 on 38 vertices. Since the numbers grow
exponentially, bridging the gap between 38 and 42 seems to be intractable. One has to abandon the idea of determining first the collection
of all combinatorial $(21_4)$ configurations and in the second step filtering out configurations that admit geometric realization.  

It seems wiser to set up a more modest goals that we state as a problem.

\begin{problem}
Determine all geometric $(21_4)$ configurations with non-trivial geometric symmetry.
\end{problem}

\begin{question}
Does there exist a geometric $(21_4)$ configuration that has no polycyclic realization?
\end{question}


\subsection{Changing the voltage group} \label{sect:voltages}


Another way to generalize $\mathrm B(21_4)$ is to change the voltage group of RLG(B) from $\mathbb Z_3$ 
to $\mathbb Z_m$, for some $m>3$. This is equivalent to saying that one expects an infinite series of configurations with rotational 
symmetry of order $m$, all with analogous structure to that of $\mathrm B(21_4)$.

We made a number of experiments for constructing such examples, using both our synthetic method in Section~\ref{sect:synthetic} 
and the procedure implemented in \emph{Mathematica} described in Section~\ref{sect:analytic}. It is clear that there likely are a number of infinite families of similar configurations. 
However, to find a proof we have to understand 
better the structure of these configurations. This is a subject of future research. Based on preliminary experiments, we conjecture

\begin{conjecture}
The following parameter values lead to geometric $(7m_{4})$ configurations, which can be realized polycyclically over $\mathbb{Z}_{m}$:
\begin{multline}\label{eq:paramSet1}
F_{1}(m; a,b) = \{a, c, d, e, f, g, q, a', c', d', e', f', g', q', t\} = \\
 \{a, -a, a, b, b, -a, -2 a, a, -a, a, b, b, -a, -2 a, 0\} \end{multline}
 for $a \geq b$,$1 \leq a, b \leq m/2$, and

\begin{equation}\label{eq:paramSet2}
F_{2}(m; a,b) = \{a, c, d, e, f, g, q, a', c', d', e', f', g', q', t\} = \{a,b,b,b,b,b,b,a,b,b,b,b,b,b, 0\} 
\end{equation}
 for $a\geq b$, $a \neq m/2$,
although there may be other constraints on $a$ and $b$ that have not yet been identified. 
\end{conjecture}
We conjecture that there are other valid as-yet-unidentified parameter families as well.

Note that the configuration $F_{1}(3; 1, 1)$, with the parameters given in  \eqref{eq:paramSet1}, is isomorphic to $B(21_{4})$. 
%
%

\begin{figure}[!h]
\begin{center}
\subcaptionbox{$F_{1}(4; 1,1)$
\label{fig:Fam1-4-1-1}}[.3\textwidth]{
\includegraphics[width = \linewidth]{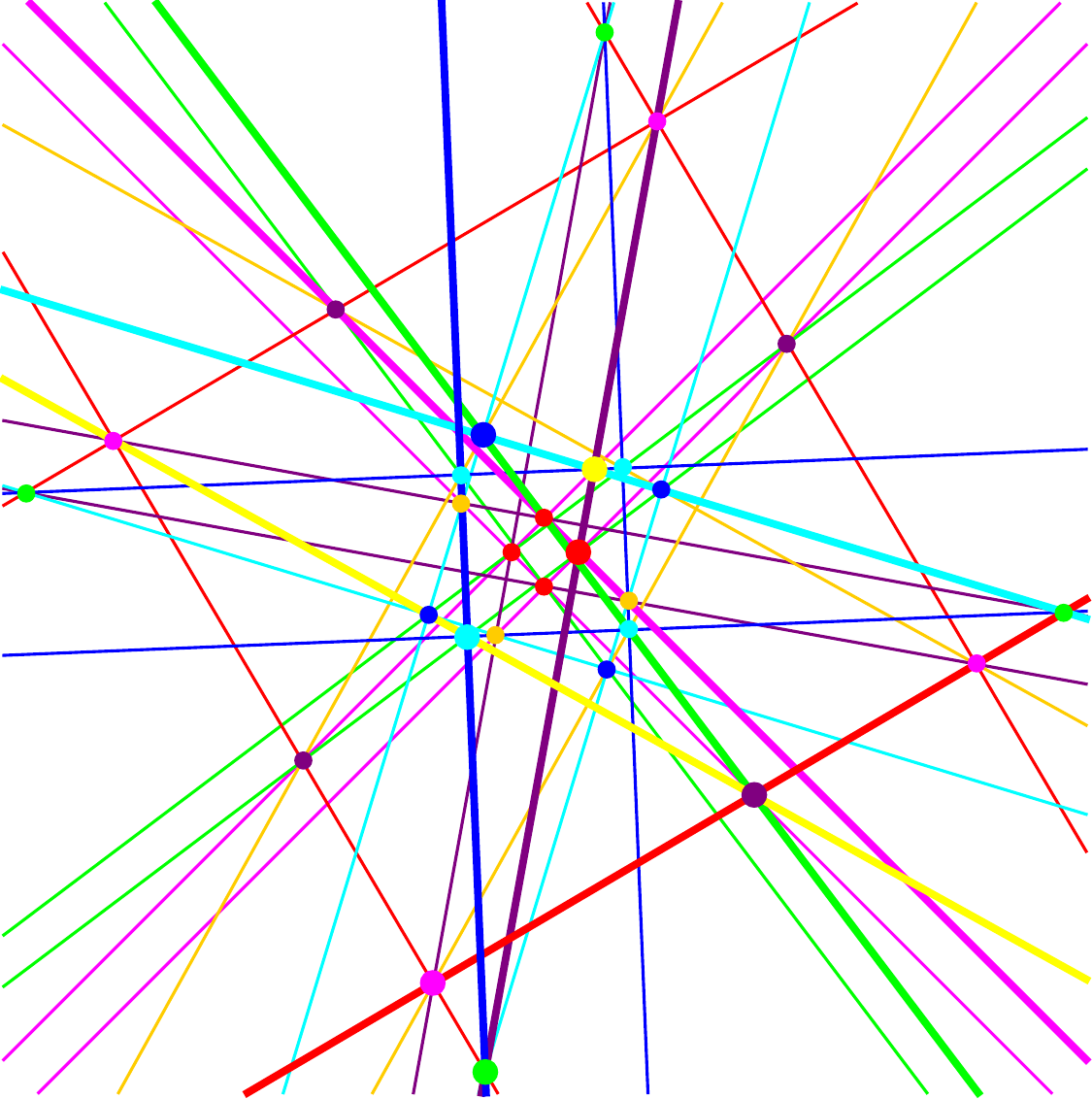}} 
\hfill
\subcaptionbox{$F_{1}(5; 2,2)$
\label{fig:Fam1-5-2-2}}[.3\textwidth]{
\includegraphics[width = \linewidth]{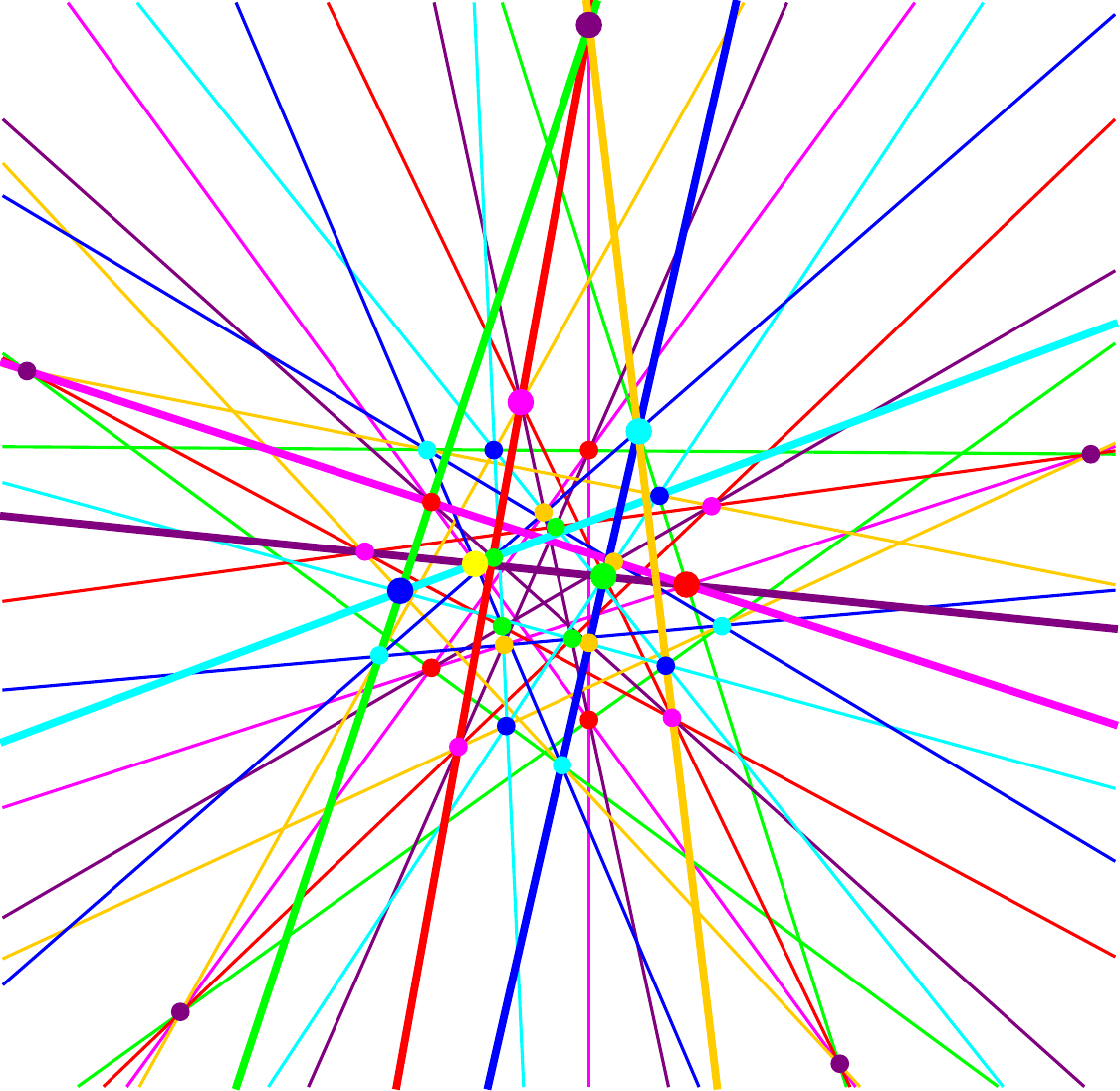}}
\hfill
\subcaptionbox{$F_{1}(6; 2,1)$
\label{fig:Fam1-4-1-1}}[.3\textwidth]{
\includegraphics[width = \linewidth]{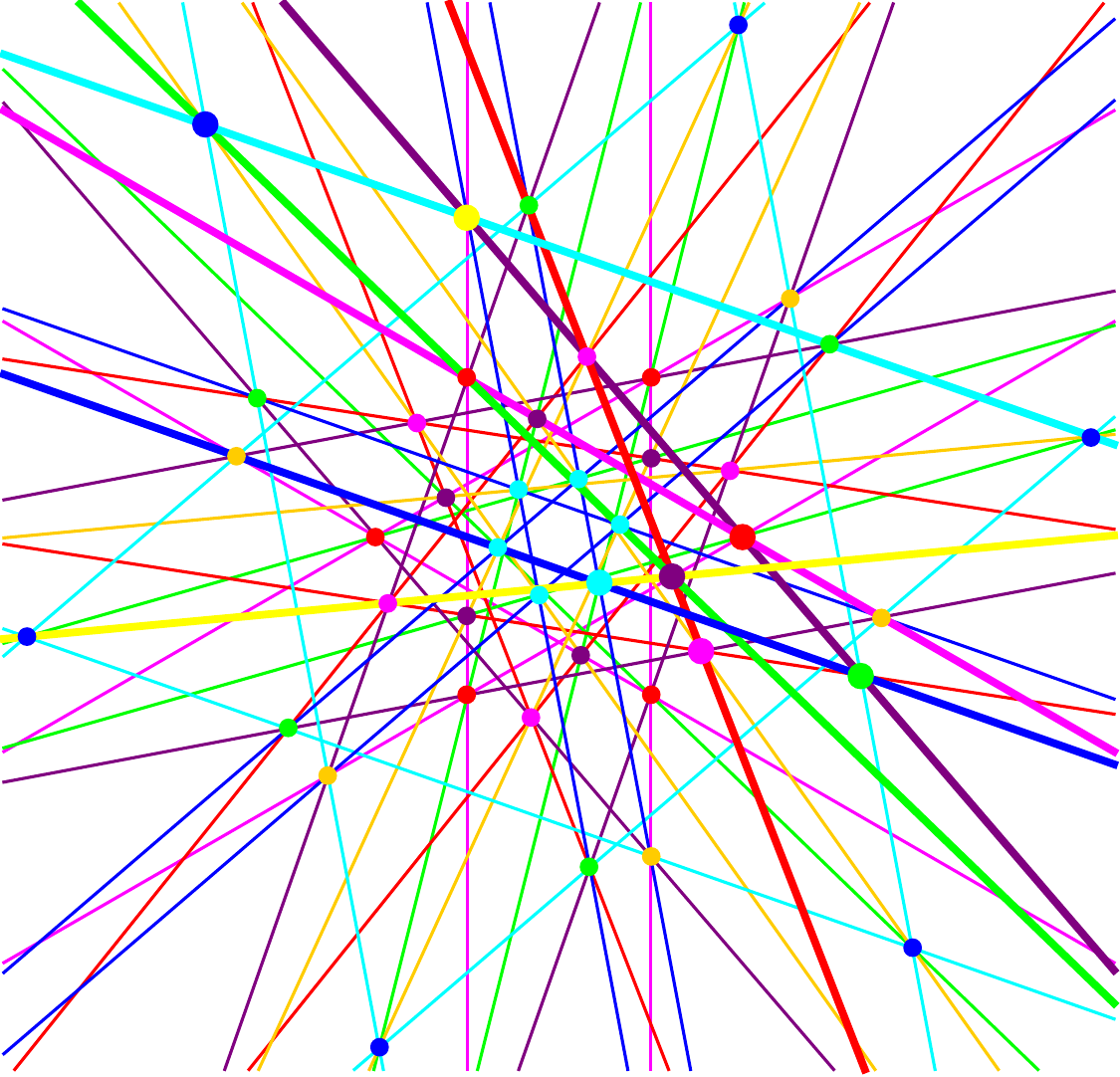}}
\caption{
Examples of configurations  in family $F_{1}(m;a,b)$.
}
\label{fig:paramFamily1}
\end{center}
\end{figure}

\begin{figure}[!h]
\begin{center}

\subcaptionbox{$F_{2}(4; 3,1)$
\label{fig:Fam3-4-3-1}}[.3\textwidth]{
\includegraphics[width = \linewidth]{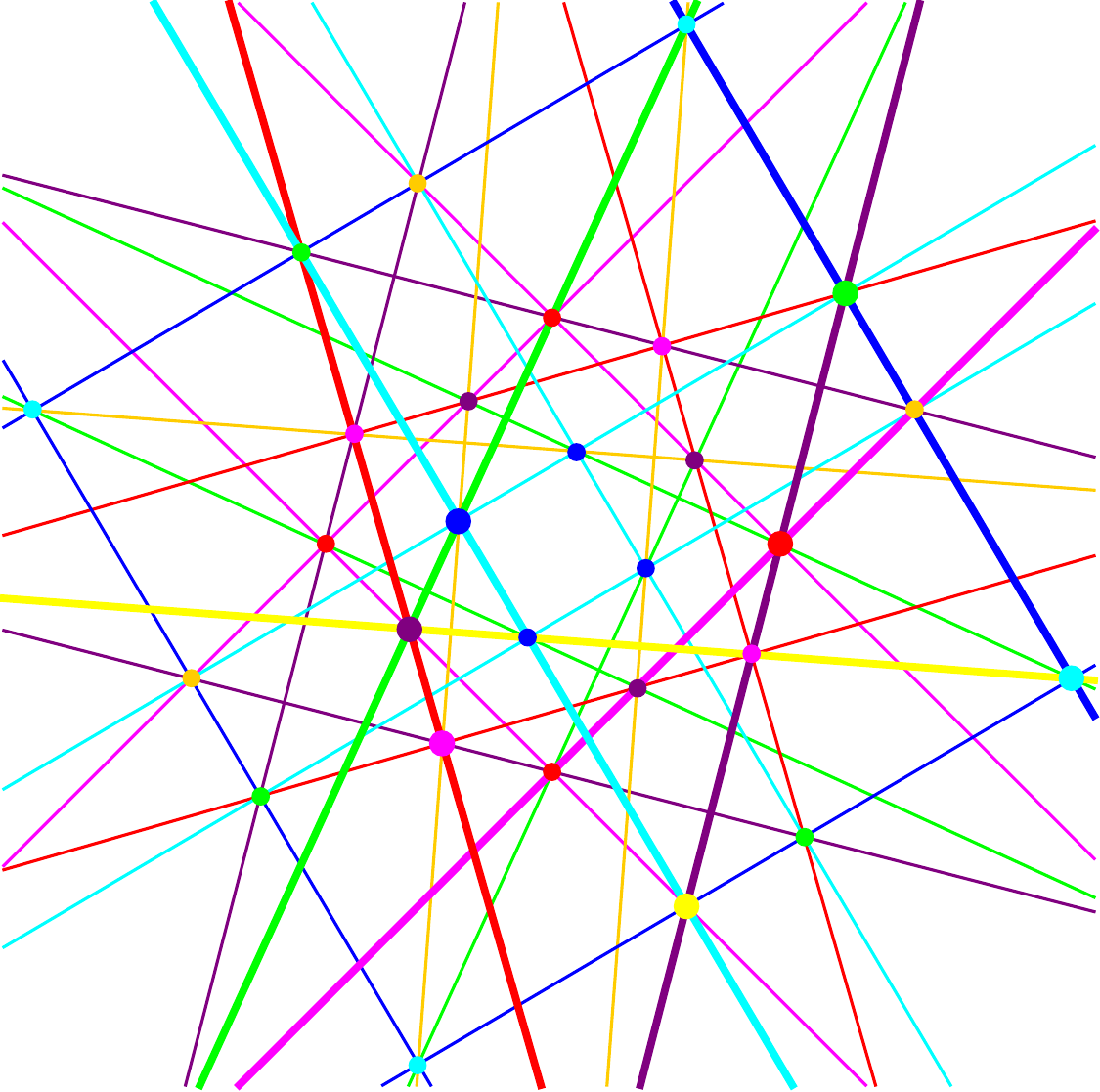}
}
\hfill 
\subcaptionbox{$F_{2}(5; 3,1)$
\label{fig:Fam3-5-3-1}}[.3\textwidth]{
\includegraphics[width = \linewidth]{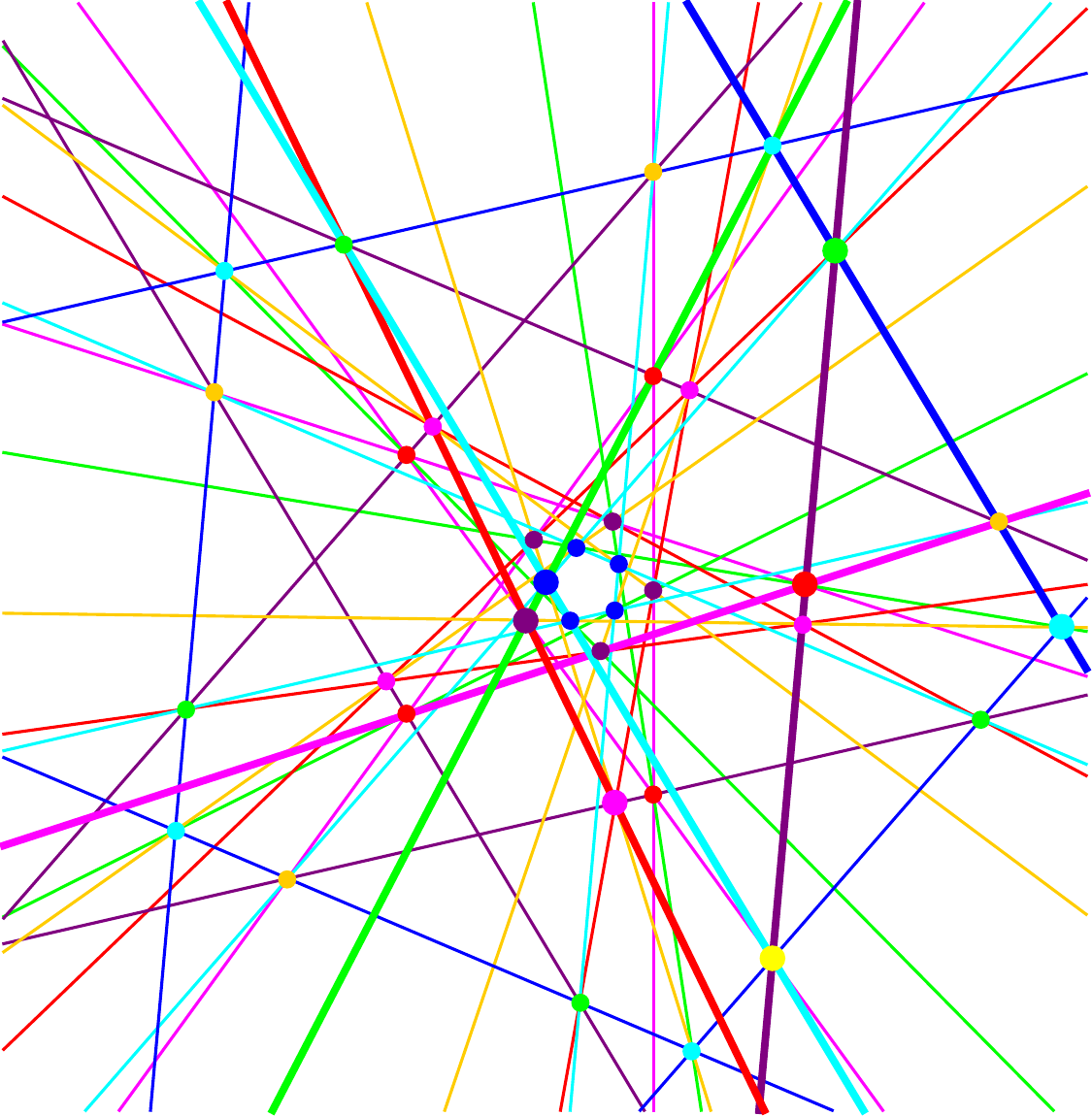}
}
\hfill
\subcaptionbox{$F_{2}(6; 4,1)$
\label{fig:Fam3-6-4-1}}[.3\textwidth]{
\includegraphics[width = \linewidth]{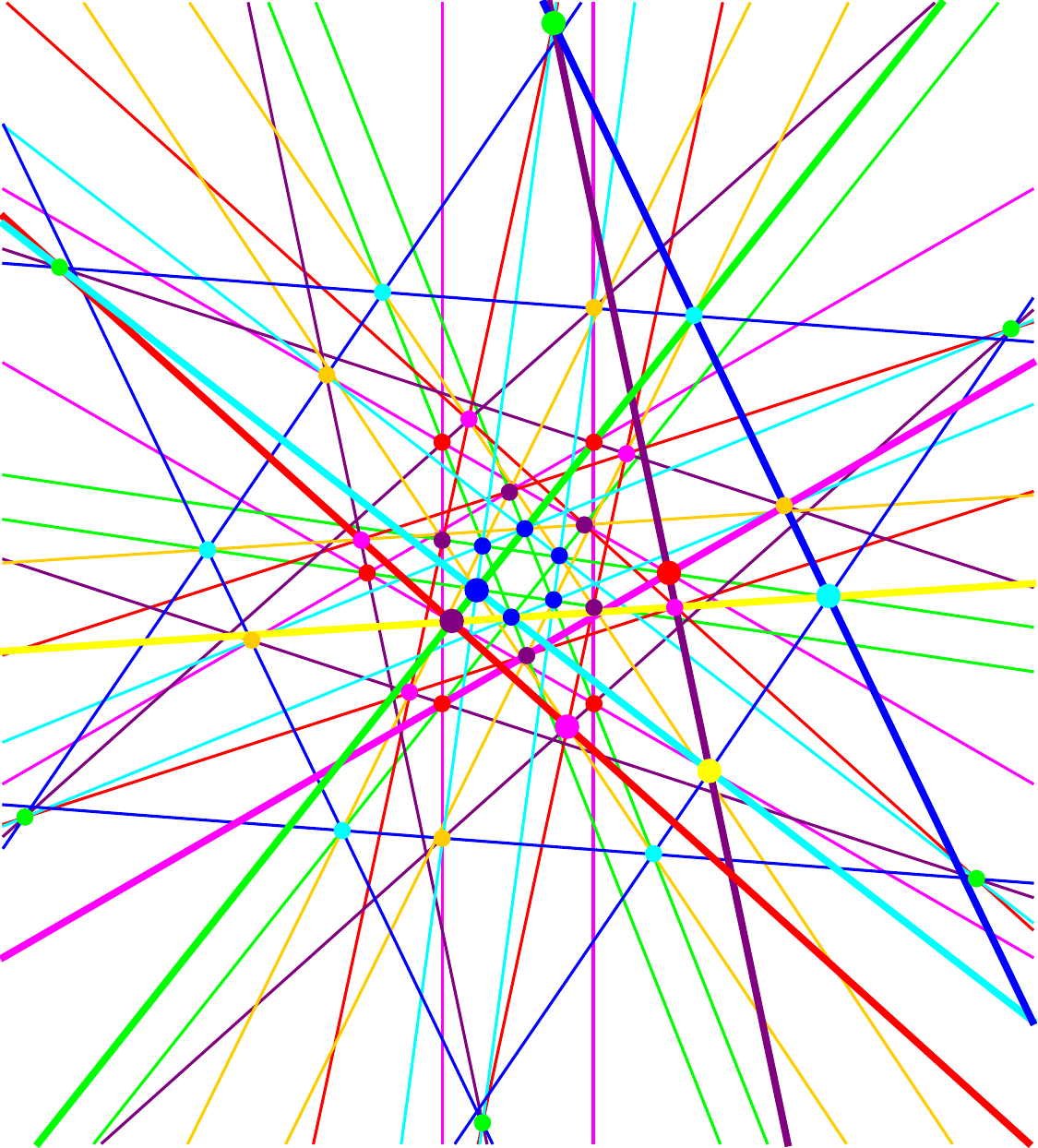}
}
\caption{
Examples of configurations  in family $F_{2}(m;a,b)$.
}
\label{fig:paramFamily2}
\end{center}
\end{figure}

We have numerically verified the existence of configurations in both parameter families for many values of $m$; they have been verified exactly for $m = 4, 6$ using the process given in section \ref{sect:analytic} (with the corresponding parameter values, naturally). Figures \ref{fig:paramFamily1} and \ref{fig:paramFamily2} show several examples of such configurations, in both conjectured parameter families. There are many other such examples.


\section*{Acknowledgements}
G\'abor G\'evay is supported by the Hungarian National Research, Development and Innovation Office, OTKA grant No.\ SNN 132625.
Toma\v{z} Pisanski is supported in part by the Slovenian Research Agency (research program P1-0294 and research projects J1-1690, N1-0140, J1-2481).

\bibliographystyle{plain}
\bibliography{b21paper}

\end{document}